\documentclass[11pt]{amsart}
\numberwithin{equation}{section}
\usepackage{amssymb}
\usepackage{amsmath}
\usepackage{latexsym}
\def\C{\mathbb C}
\def\N{\mathbb N}

\newtheorem{la}{Lemma}[section]
\newtheorem{thm}{Theorem}[section]
\newtheorem{que}{Question}[section]
\newtheorem{cor}{Corollary}[section]
\theoremstyle{definition}
\newtheorem{defin}{Definition}[section]
\theoremstyle{remark}

\title[Hyperbolic Meromorphic Functions ]{Dynamics of Hyperbolic Meromorphic Functions}
\author{Zheng Jian-Hua}
\address{Department of Mathematical Sciences, Tsinghua University, P. R. China}
\email{jzheng@math.tsinghua.edu.cn}
\begin{document}
\begin{abstract}
A definition of hyperbolic meromorphic functions is given and then
we discuss the dynamical behavior and the thermodynamic formalism of
hyperbolic functions on the Julia set. We prove the important
expanding properties for hyperbolic functions on the complex plane
and with respect to the euclidean metric. We establish the Bowen
formula for hyperbolic functions on the complex plane, that is, the
Poincare exponent equals to the Hausdorff dimension of the radial
Julia set and furthermore, we prove that all the results in the
Walters' theory hold for hyperbolic functions on the Riemann sphere.
\end{abstract}
\maketitle
\section{Introduction}

Let $f(z)$ be a meromorphic function which is transcendental or
rational with degree at least two. Consider the $n$-th iterate of
$f(z)$ defined by $$f^0(z)=z,\ f^n(z)=f(f^{n-1}(z))=f^{n-1}(f(z)).$$
It is clear that if $f(z)$ is transcendental, that is, $\infty$ is
an essential singular point of $f(z)$, then $f^n(z)$ is meromorphic
only on $\C\setminus \cup_{j=1}^{n-1}f^{-j}(\infty)$.

Let $\mathcal{F}_f$ be the Fatou set of $f(z)$ defined by
$$\mathcal{F}_f=\{z\in\widehat{\mathbb{C}}: \{f^n\}\ {\rm is\ well\
defined\ and\ normal\ at\ a\ neighborhood\ of}\ z\}$$ and
$\widehat{\mathcal{J}}_f$ the Julia set of $f(z)$, that is,
$\widehat{\mathcal{J}}_f=\widehat{\mathbb{C}}\setminus
\mathcal{F}_f$ and set
$\mathcal{J}_f=\widehat{\mathcal{J}}_f\setminus\{\infty\}$. If
$f(z)$ is rational, then it is possible that $\infty$ is in
$\mathcal{F}_f$ and in this case
$\mathcal{J}_f=\widehat{\mathcal{J}}_f$; If $f(z)$ is
transcendental, then $\infty$ must be in $\widehat{\mathcal{J}}_f$.
The prepoles of $f(z)$ are important points which we have to take
more care of in this note. A point $z_0$ in $\mathbb{C}$ is called
prepole of $f(z)$ if for some $n\geq 1$, $f^n(z_0)=\infty$. Set
$$\mathcal{J}_f(\infty)=\bigcup_{j=0}^\infty f^{-j}(\infty).$$ If
$\infty$ is not a Picard exceptional value of $f(z)$, then
$\widehat{\mathcal{J}}_f=\overline{\mathcal{J}_f(\infty)}.$ However,
for a transcendental entire function $f(z)$, we take more care of
$\mathcal{J}_f$ instead of $\widehat{\mathcal{J}}_f$.

By $\widehat{{\rm sing}}(f^{-1})$ we denote the closure of the set
of critical and asymptotic values of $f(z)$ (including $\infty$ in
our consideration). If $f(z)$ has a pole with multiplicity at least
two, then $\infty$ is a critical value of $f(z)$; A rational
function has no asymptotic value; A transcendental entire function
has $\infty$ as an asymptotic value of it. By ${\rm sing}(f^{-1})$
we denote the closure in $\mathbb{C}$ of the set of finite critical
and asymptotic values of $f(z)$ (i.e., excluding $\infty$ from our
consideration) and by $\mathcal{P}(f)$ the post-singular set defined
to be the closure in $\widehat{\mathbb{C}}$ of
$$\bigcup_{n=0}^\infty f^n({\rm
sing}(f^{-1})\setminus \cup_{j=0}^{n-1}f^{-j}(\infty))$$ and set
$\widehat{\mathcal{P}}(f)=\mathcal{P}(f)\cup\widehat{{\rm
sing}}(f^{-1})$. One has proved that $\mathcal{P}(f)$ plays an
important role in study of dynamical aspects of meromorphic
functions. Some of the known results which describe the fact will be
listed in the next section.

In this note, one of our main purposes is to discuss the dynamical
behavior of hyperbolic meromorphic functions. Let us begin with
definition of hyperbolic meromorphic functions. The definition of
hyperbolic rational functions is clear. However, the transcendental
case is not so and two different definitions have been given out in
other papers. A transcendental meromorphic function with
\begin{equation}\label{1.1}\mathcal{P}(f)\cap\widehat{\mathcal{J}}_f=\emptyset\end{equation} is defined
by Rippon and Stallard \cite{RipponStallard} to be hyperbolic and
with
\begin{equation}\label{1.1+}{\rm dist}(\mathcal{P}(f), \mathcal{J}_f)>0\end{equation}
by Mayer and Urb\'anski \cite{MayerUrbanski} to be topologically
hyperbolic and to be hyperbolic if, in addition, for some $c>0$ and
some $\lambda>1$,
\begin{equation}\label{1.2}|(f^n)'(z)|\geq c\lambda^n,\ \forall\
n\geq 1,\ \forall\ z\in
\mathcal{J}_f\setminus\mathcal{J}_f(\infty).\end{equation} Since
$\infty\in\widehat{\mathcal{J}}_f$, the condition (\ref{1.1})
implies that $\infty\not\in \mathcal{P}(f)$, that is,
$\mathcal{P}(f)$ is bounded and then it satisfies (\ref{1.1+}). The
condition (\ref{1.1+}) does not imply (\ref{1.1}), but it together
with ${\rm sing}(f^{-1})$ being bounded yields (\ref{1.1}) (see
\cite{RipponStallard}).

We put special attention on $\infty$ and in order to avoid
occurrence of confusion about definitions of hyperbolic meromorphic
functions, we give out the following definition.

\begin{defin}\label{def1.1}\ Let $f(z)$ be a meromorphic function in $\mathbb{C}$ transcendental or rational with degree
at least two.

(1) $f(z)$ is called hyperbolic on the Riemann sphere (or with
respect to the sphere metric in order to emphasize the considered
metric) if
\begin{equation}\label{1.3}\widehat{\mathcal{P}}(f)\cap\widehat{
\mathcal{J}}_f=\emptyset;\end{equation}

(2) $f(z)$ is called hyperbolic on the complex plane if (\ref{1.1})
holds;

(3) $f(z)$ is called hyperbolic with respect to the Euclidean metric
if it satisfies (\ref{1.1+}) and (\ref{1.2}).
\end{defin}

The reasons for these names can be found from Theorem \ref{thm2.2},
Theorem \ref{thm2.3} and Theorem \ref{thm2.3+}. A hyperbolic
function on the complex plane may not be hyperbolic with respect to
the Euclidean metric; A function is hyperbolic with respect to the
Euclidean metric if and only if it is expanding on
$\mathcal{J}_f\setminus\mathcal{J}_f(\infty)$ with respect to the
Euclidean metric (see Theorem \ref{thm2.3+}); A hyperbolic function
with respect to the Euclidean metric may not be hyperbolic with
respect to the sphere metric. A function is hyperbolic on the
Riemann sphere if and only if it is expanding on
$\mathcal{J}_f\setminus\mathcal{J}_f(\infty)$ with respect to the
sphere metric (see Theorem \ref{thm2.3}).

For the sake of simplicity, let us introduce following notations: we
denote by $\mathcal{H}(\widehat{\mathbb{C}})$ the set of all
hyperbolic meromorphic functions on the Riemann sphere; by
$\mathcal{H}(\mathbb{C})$ the set of all hyperbolic meromorphic
functions on the complex plane; by $\mathcal{H}(Eu)$ the set of all
hyperbolic meromorphic functions with respect to the Euclidean
metric. Then it is obvious that
$\mathcal{H}(\widehat{\mathbb{C}})\subset\mathcal{H}(\mathbb{C})$
and $\mathcal{H}(Eu)\cap\mathcal{B}\subset\mathcal{H}(\mathbb{C}),$
where $\mathcal{B}$ is the set of all bounded-type functions, that
is, functions with bounded ${\rm sing}(f^{-1})$. However, we do not
know if $\mathcal{H}(\widehat{\mathbb{C}})\subset \mathcal{H}(Eu).$

The definition of hyperbolic meromorphic functions on the Riemann
sphere coincides with that of hyperbolic rational functions, but the
others are deferent from that. Notice that if $f(z)$ is
transcendental, $\infty$ must be in $\widehat{ \mathcal{J}}_f$ and
thus $\infty$ is not a singular value of a transcendental hyperbolic
meromorphic function in $\mathcal{H}(\widehat{\mathbb{C}})$, that
is, the inverse of $f(z)$ has no singularities over $\infty$.
Therefore, for $f\in\mathcal{H}(\widehat{\mathbb{C}})$ and for all
sufficiently large $R=R(f)$,
$$f^{-1}(\mathcal{D}_R)=\bigcup\limits_{j=1}^\infty V_j$$ where
$\mathcal{D}_R=\{|z|>R\}\cup\{\infty\}$ such that each $V_j$ is a
simply connected domain containing a pole of $f(z)$ and
$f:V_j\rightarrow \mathcal{D}_R$ is univalent. Furthermore, each
$V_j$ is bounded and surrounded by an analytic Jordan curve.

We discuss the possibility of that $\infty$ is an asymptotic value.
A classical theorem of Iversen \cite{Iversen} implies that a
transcendental meromorphic function with only finitely many poles
has $\infty$ as an asymptotic value and hence such functions are not
hyperbolic on the Riemann sphere. Furthermore, if $f(z)$ is
transcendental and in $\mathcal{H}(\widehat{\mathbb{C}})$, then it
is of bounded type, that is, ${\rm sing}(f^{-1})$ is bounded, and
$\infty$ is not an asymptotic value, and thus in view of a result of
Teichm\"uller \cite{Teichmuller} (cf. Proposition 7.1 of Bergweiler,
Rippon and Stallard \cite{BerRiSt}), $f(z)$ has zero Nevanlinna
deficiency at $\infty$, that is, $\delta(\infty,f)=0$. Consequently,
we have shown that a transcendental hyperbolic meromorphic function
in $\mathcal{H}(\widehat{\mathbb{C}})$ has $\delta(\infty,f)=0$, and
unforturely  a meromorphic function with $\delta(\infty,f)>0$ in
$\mathcal{H}(\mathbb{C})$ is not in
$\mathcal{H}(\widehat{\mathbb{C}})$. However, a transcendental
meromorphic function with the Nevanlinna deficiency
$\delta(\infty,f)=1$ may not have $\infty$ as its asymptotic value,
which was proved in Hayman \cite{Hayman} and Ter-Israelyan
\cite{Ter-Israelyan}, and then such functions are not of bounded
type.

Let us observe examples of transcendental hyperbolic meromorphic
functions. For $0<\lambda<1$, $\lambda\tan z$ is hyperbolic in
$\mathcal{H}(\widehat{\mathbb{C}})$ and for sufficiently small
$\lambda>0$, $\lambda\tan(\pi\sin z)$ is also hyperbolic in
$\mathcal{H}(\widehat{\mathbb{C}})$ and of infinite order. The
function
\begin{equation}\label{1.5}f_{p,\lambda}(z)=\lambda\sum_{n=p^2}^\infty\left(\frac{1}{n^p-z}-\frac{1}{n^p+z}\right),\
p\in\N\end{equation} is hyperbolic in
$\mathcal{H}(\widehat{\mathbb{C}})$ for
$0<\lambda<\frac{p^{4p-1}}{10^4\log p}$ and $p\geq 6$, which can be
proved from Stallard \cite{Stallard}. In general, we can consider
the function
$$\frac{f'(z)}{f(z)},\ {\rm where}\ f(z)=\prod_{n=1}^\infty
\left(1-\frac{z}{a_n}\right)\ {\rm and}\ \sum_{n=1}^\infty
\frac{1}{|a_n|}<\infty.$$ Many examples in $\mathcal{H}(Eu)$ can be
found in Chapter 3 of Mayer and Urb\'anski \cite{MayerUrbanski}.

The main purpose of this paper has two ones: one is to study the
topological dynamics of hyperbolic meromorphic functions and the
other is to study the thermodynamical formalism of hyperbolic
meromorphic functions. We have understood relatively clearly these
two aspects of hyperbolic rational functions and many excellent
results have been revealed. This motivates us to investigate the
transcendental hyperbolic case. We shall find that almost of the
dynamical behaviors of transcendental hyperbolic meromorphic
functions in $\mathcal{H}(\widehat{\mathbb{C}})$ are the same as
those of hyperbolic rational functions. For example, a Cantor Julia
set can be expressed by symbolic shift automorphism and, among other
things, we shall discuss some relationship between the topological
behavior of Julia set and $\infty$ being not an asymptotic value.
The expanding properties are an important characteristic of
hyperbolic functions. We shall prove that the function in
$\mathcal{H}(\widehat{\mathbb{C}})$ and in $\mathcal{H}(Eu)$ are
expanding respectively with respect to the sphere metric and the
euclidean metric. We shall prove the existence for suitable $t>0$
and the strictly decreasing property in $t$ when it exists of the
pressure function $P(f,t)$ of hyperbolic functions in
$\mathcal{H}(\mathbb{C})$, and furthermore, we shall deduce the
Bowen formula, that is, the Poincare exponent equals to the
Hausdorff dimension of the radial Julia set (see Theorem 3.3 below),
which refines precisely a result of Stallard \cite{Stallard} (see
Theorem 3.2 below) and by noting that
$\mathcal{H}(\widehat{\mathbb{C}})\subset\mathcal{H}(\mathbb{C})$,
we shall get that for a function in
$\mathcal{H}(\widehat{\mathbb{C}})$, the Poincare exponent exactly
equals to the Hausdorff dimension of the Julia set. Kotus and
Urbanski \cite{KotusUrbanski} used the Walters' theory
\cite{Walters} to prove the Bowen formula for the functions in
$\mathcal{H}(\widehat{\mathbb{C}})$ under the assumption of that the
functions are so called strongly regular. The final part of this
paper is devoted to proving that every function in
$\mathcal{H}(\widehat{\mathbb{C}})$ satisfies each item of
conditions in Walters's theory except the expanding condition, while
in terms of Theorem 2.2, for some $N$, $f^N$ satisfies the expanding
condition. Thus we prove that all the results in Walters' theory
\cite{Walters} hold for every function in
$\mathcal{H}(\widehat{\mathbb{C}})$.

 Finally, we introduce
basic notations which will be often used in the paper. By
$d_\infty(a,b)$ we denote the sphere distance between two points $a$
and $b$ on $\widehat{\mathbb{C}}$ and by $d(a,b)$ the Euclidean
distance between two points $a$ and $b$ on $\mathbb{C}$. We mean by
$B_\infty(a,r)$ and $B(a,r)$ the disks centered at $a$ with radius
$r$ under respectively the sphere metric $d_\infty$ on
$\widehat{\mathbb{C}}$ and the Euclidean metric  $d$ on
$\mathbb{C}$. Let $f(z)$ be a meromorphic function on a subdomain
$D$ of $\widehat{\mathbb{C}}$. By $f^\times(z)$ we denote the
derivative of $f(z)$ with respect to the sphere metric on $D$, that
is, for $z\in D\setminus\{\infty\}$,
$$f^\times(z)=\frac{|f'(z)|(1+|z|^2)}{1+|f(z)|^2},\
f(z)\not=\infty;$$
$$f^\times(z)=\lim_{\zeta\rightarrow z}f^\times(\zeta),\ f(z)=\infty$$
and for $z=\infty$,
$$f^\times(\infty)=\lim_{\zeta\rightarrow 0}F^\times(\zeta),\
F(\zeta)=1/f(1/\zeta).$$ For example, consider $f(z)=\frac{z^2}{
(z-1)(z-2)}$ and $g(z)=2z+\frac{1}{z}$ and then
$f^\times(\infty)=((1-\zeta)(1-2\zeta))^\times(0)=\frac{3}{2}$ and
$g^\times(\infty)=\left(\frac{\zeta}{2+\zeta^2}\right)^\times(0)=\frac{1}{2}.$
For a transcendental meromorphic function, we cannot consider its
derivative at $\infty$ with respect to the sphere metric.

A point $z_0\in\widehat{\mathbb{C}}$ is called a periodic point of
$f(z)$ with period $p$ if $f^p(z_0)=z_0$ and $f^j(z_0)\not=z_0$ for
$1\leq j\leq p-1$. At periodic points with period $p$,
$(f^p)^\times(z)=|(f^p)'(z)|$. A periodic point $z_0$ is called in
turn attracting, indifferent and repelling if $(f^p)^\times(z_0)$ is
less than, equals or greater than $1$.

\vskip 1cm
\section{Expanding and Topological Dynamics of Hyperbolic Functions}

One has revealed the close connections between the components of the
Fatou set and singular values of meromorphic functions. For our
purpose of this paper, we recall some of them. The components of the
Fatou set are classified into (pre)periodic domains and wandering
domains. A component $U$ of the Fatou set $\mathcal{F}_f$ of $f$ is
called wandering if $f^n(U)\cap f^m(U)=\emptyset$ for $n\not=m$;
periodic if for some $n>0$, $f^n(U)\cap U\not=\emptyset$ and in this
case, in fact we have $f^n(U)\subseteq U$ and the smallest $n$ such
that the inclusion holds is called period of $U$; pre-periodic if
for some $n>0$, the component containing $f^n(U)$ is periodic, but
$U$ is not periodic. The periodic domains are classified into
attracting domains, parabolic domains, Siegel disks, Herman rings
and Baker domains. Every cycle of attracting domains and parabolic
domains contains at least one singular value; the boundaries of
Siegel disks and Herman rings are subset of the post-singular set
$\widehat{\mathcal{P}}(f)$. About the cases of the Baker domains and
wandering domains, we have the following

\begin{thm}\label{thm2.1}\ \ Let $f(z)$ be a transcendental
meromorphic function.  Then the following statements hold:

(1) if $U$ is a Baker domain of $f(z)$ with period $p$, then
$a=\lim\limits_{n\rightarrow\infty}f^{pn}|_U$ is in the derived set
of $\mathcal{P}_p(f)\cap \mathcal{J}_f$ on $\widehat{\mathbb{C}}$,
where
$$\mathcal{P}_p(f)=\bigcup_{j=0}^{p-1}f^j({\rm
sing}(f^{-1})\setminus \cup_{i=0}^{j-1}f^{-i}(\infty));$$

(2) if $U$ is a wandering domain, then all limit values of $\{f^n\}$
on $U$ are in the derived set of
$\bigcup_{p=1}^\infty\mathcal{P}_p(f)\cap \mathcal{J}_f$ on
$\widehat{\mathbb{C}}$.\end{thm}

Theorem \ref{thm2.1} was proved in Zheng\cite{Zheng} and the result
(2) for $f(z)$ being an entire function was verified by Bergweiler
et al \cite{Bergweileretal}.

Now we establish equivalent results of hyperbolic meromorphic
functions which are well-known for the case of hyperbolic rational
functions.

\begin{thm}\label{thm2.2}\ \ Let $f(z)$ be a
meromorphic function transcendental or rational with degree at least
two. Then the following statements are equivalent:

(1) $f(z)$ is hyperbolic on the Riemann sphere, that is,
$f\in\mathcal{H}(\widehat{\mathbb{C}})$;

(2) each point of $\widehat{{\rm sing}}(f^{-1})$ is attracted to a
cycle of attracting periodic points and $f(z)$ has only finitely
many attracting periodic points;

(3) $f(z)$ is expanding on $\widehat{\mathcal{J}}_f$ with respect to
the sphere metric, that is, for some $m$, some $\lambda>1$ and some
$\delta>0$, we have
\begin{equation}\label{2.1}d_\infty(f^m(z),f^m(w))\geq \lambda d_\infty(z,w)\end{equation} whenever
$z$ and $w$ are in a common component of
$f^{-m}(B_\infty(a,\delta))$ for some $a\in
\widehat{\mathcal{J}}_f$.
%and $d_\infty(z,w)<\delta$.
\end{thm}

{\bf Proof.}\ \ It suffices to prove our results for the case of
transcendental meromorphic functions.

To prove (2) from (1). Assume that $f(z)$ is hyperbolic on the
Riemann sphere and then $\widehat{{\rm sing}}(f^{-1})={\rm
sing}(f^{-1})$ and ${\rm sing}(f^{-1})$ is bounded and furthermore
${\rm sing}(f^{-1})\subset\mathcal{F}_f$ and $f(z)$ has no parabolic
domains, Siegel disks and Herman rings. In view of Theorem
\ref{thm2.1}, $f(z)$ has no Baker domains and wandering domains.
This yields that each point of ${\rm sing}(f^{-1})$ is attracted to
a cycle of attracting periodic points. Suppose that $f(z)$ has
infinitely many attracting periodic points and we take a sequence
$\{a_n\}$ of such points. It is easily seen that all limit points of
$\{a_n\}$ are in $\widehat{\mathcal{J}}_f$. Since $\{a_n\}\subset
\mathcal{P}(f)$, all limit points of $\{a_n\}$ are in
$\widehat{\mathcal{J}}_f\cap \widehat{\mathcal{P}}(f)$ so that
$\widehat{\mathcal{J}}_f\cap
\widehat{\mathcal{P}}(f)\not=\emptyset$, which contradicts the
condition (1). Thus we have proved that $f(z)$ has only finitely
many attracting periodic points. Therefore the implication (1)
$\Rightarrow$ (2) is completed.

To proceed the proof of (2) $\Rightarrow $ (1). Since $\infty$ must
be in $\widehat{\mathcal{J}}_f$, the condition (2) implies that
$\infty$ is not in $\widehat{{\rm sing}}(f^{-1})$ and hence,
$\widehat{{\rm sing}}(f^{-1})={\rm sing}(f^{-1})$ and ${\rm
sing}(f^{-1})$ is bounded. Since every limit point of ${\rm
sing}(f^{-1})$ is still in ${\rm sing}(f^{-1})$, the sphere distance
from ${\rm sing}(f^{-1})$ to $\widehat{\mathcal{J}}_f$ is positive
and furthermore $\widehat{\mathcal{P}}(f)$ is compact on
$\mathbb{C}$. Thus (1) easily follows from (2).

To prove (3) from (1). Set
$$\eta=\frac{1}{4}d_\infty(\widehat{\mathcal{J}}_f,\widehat{\mathcal{P}}(f))>0.$$
We can take finitely many points $a_j (1\leq j\leq d)$ on
$\widehat{\mathcal{J}}_f$ such that $B_\infty(a_j,\eta)\ (1\leq
j\leq d)$ form a covering of $\widehat{\mathcal{J}}_f$ and then
$f^{-1}(B_\infty(a_j,\eta))\ (1\leq j\leq d)$ form a covering of
$\mathcal{J}_f$. $f^{-1}(B_\infty(a_j,4\eta))\ (1\leq j\leq d)$ do
not contain $\infty$ and points of $\mathcal{P}(f)$. Every component
of $f^{-1}(B_\infty(a_j,4\eta))\ (1\leq j\leq d)$ is simply
connected and on it, for every $n$, $f^{-n}(\zeta)$ consists of
single-valued analytic branches, denoted by $f^{-n}_{j,k}(\zeta)$ on
$B_\infty(a_j,4\eta)\ (1\leq j\leq d; 1\leq k< \infty)$. Thus for
every $n$, $f_{j,k}^{-n}(B_\infty(a_j,4\eta))\ (1\leq j\leq d; 1\leq
k< \infty)$ is a covering of $\mathcal{J}_f\setminus
\mathcal{J}_f(\infty)$. We claim that
$(f^{-n}_{j,k})^\times(\zeta)\rightarrow 0$ as $n\rightarrow\infty$
uniformly on $B_\infty(a_j,\eta)$. Suppose on contrary the claim
fails. Then there exist two sequences of positive integers $\{n_m\}$
and $\{k_m\}$ with $n_m\rightarrow\infty$ as $m\rightarrow\infty$
and a sequence of complex numbers $\{\zeta_m\}$ in
$B_\infty(a_j,\eta)$ such that
$(f_{j,k_m}^{-n_m})^\times(\zeta_m)>\varepsilon_0>0$ and further,
$(f_{j,k_m}^{-n_m})^\times(a_j)>c\varepsilon_0>0$ for an absolute
constant $c$. In view of the Koebe Theorem, we have
$$B_\infty(f_{j,k_m}^{-n_m}(a_j),\frac{1}{4}c\varepsilon_0\eta)\subset
B_\infty(f_{j,k_m}^{-n_m}(a_j),\frac{1}{4}(f_{j,k_m}^{-n_m})^\times(\zeta_m)\eta)\subset
(f_{j,k_m}^{-n_m})(B_\infty(a_j,\eta))$$ so that
$$f^{n_m}(B_\infty(f_{j,k_m}^{-n_m}(a_j),\frac{1}{4}c\varepsilon_0\eta))\subset
B_\infty(a_j,\eta).$$ We assume without any loss of generalities
that $f_{j,k_m}^{-n_m}(a_j)\rightarrow b_j\in
\widehat{\mathcal{J}}_f$ and then for all sufficiently large $m$,
$B_\infty(b,\frac{1}{8}c\varepsilon_0\eta)\subset
B_\infty(f_{j,k_m}^{-n_m}(a_j),\frac{1}{4}c\varepsilon_0\eta)$ and
$f^{n_m}(b,\frac{1}{8}c\varepsilon_0\eta)\subset
B_\infty(a_j,\eta).$ This is impossible. We have proved the claim.

It is clear from the claim that for some large fixed $N$,
$(f^{-N}_{j,k})^\times(\zeta)<\frac{1}{2}$ uniformly on
$B_\infty(a_j,\eta)$ for each $j$ and each $k$. Set
$$\widehat{U}=\bigcup_{j,k}f^{-N}_{j,k}(B_\infty(a_j,\eta))$$
and $\widehat{U}$ is a neighborhood of $\mathcal{J}_f\setminus
\mathcal{J}_f(\infty)$. For an arbitrary point $z\in\widehat{U}$,
there exists a point $\zeta\in B_\infty(a_{j_0},\eta)$ such that for
some $k_0$, $f^{-N}_{j_0,k_0}(\zeta)=z$ and thus
$$(f^N)^\times(z)=\frac{1}{(f^{-N}_{j_0,k_0})^\times(\zeta)}>2.$$

Let $2\delta>0$ be the Lebesgue number of $B_\infty(a_j,\eta)\
(1\leq j\leq d)$ and then for each $a\in\widehat{\mathcal{J}}_f$, we
have a point $b\in\mathcal{J}_f\setminus \mathcal{J}_f(\infty)$ such
that $B_\infty(a,\delta)\subset B_\infty(b,2\delta)\subset
B_\infty(a_{j_0},\eta)$ for some $j_0$. For two arbitrary points $z$
and $w$ in a component $V$ of $f^{-N}(B_\infty(a,\delta))$, we have
$$d_\infty(f^N(z),f^N(w))=\int_\Gamma\frac{|d\zeta|}{1+|\zeta|^2}=\int_\gamma(f^N)^\times(\zeta)
\frac{|d\zeta|}{1+|\zeta|^2}>2d_\infty(z,w),$$ where $\Gamma$ is a
straight segment connecting $f^N(z)$ and $f^N(w)$ in
$B_\infty(a,\delta)$ and $\gamma$ is the curve connecting $z$ and
$w$ in the component $V$ such that $\Gamma=f^N(\gamma)$. We have
gotten the desired result (3).

To complete the implication (3) $\Rightarrow$ (1). Under the
condition (3), $f^{-m}$ has no singularities over
$B_\infty(a,\delta)$ for each $a\in\widehat{\mathcal{J}}_f$ and
every component of $f^{-m}(B_\infty(a,\delta))$ will lie in a disk
$B_\infty(b,\delta)$ for some $b\in\widehat{\mathcal{J}}_f$. This
implies that every $f^{-nm}$ and so every $f^{-n}$ has no
singularities over $B_\infty(a,\delta)$ for each
$a\in\widehat{\mathcal{J}}_f$ and therefore
$d_\infty(\widehat{\mathcal{J}}_f,\widehat{\mathcal{P}}(f))\geq\delta>0$.
We attain (1). \qed

The following is a basic result which comes from the unique theorem
of analytic function.

\begin{la}\label{lem2.1}\ \ Let $f(z)$ be a meromorphic function on
$\mathbb{C}$. If $\infty$ is not its singular value and for some
$R_0>0$, every component of $f^{-1}(\mathcal{D}_{R_0})$ is simply
connected and $f(z)$ is univalent on every component of
$f^{-1}(\mathcal{D}_{R_0})$, then for $r>R_0$, all but at most
finitely many components of $f^{-1}(\mathcal{D}_r)$ lie in
$\mathcal{D}_R$ for arbitrary $R>0$.\end{la}

{\bf Proof.}\ \ Suppose that there exists a positive number $R$ such
that $f^{-1}(\mathcal{D}_r)$ has infinitely many components $U_n$
which are not in $\mathcal{D}_R$. Choose a point $a$ in
$\mathcal{D}_r$ such that $f^{-1}(a)$ is infinite and then each
$U_n$ contains exactly one point $a_n$ of $f^{-1}(a)$. Obviously,
$a_n\rightarrow \infty$ as $n\rightarrow\infty$ and $U_n\cap
B(0,R)\not=\emptyset$. Hence ${\rm diam}(U_n)\rightarrow\infty$ as
$n\rightarrow\infty$. By $L$ we denote the limit set of boundaries
of $U_n$. It is easily seen that $L$ must contain an unbounded
component $\Gamma$ and $f(z)$ maps $\Gamma$ onto the circle
$\{|z|=r\}$. However, for $r>R_0$, every branch of $f^{-1}$ is
analytic on $\{|z|=r\}$ and hence each image of $\{|z|=r\}$ is
bounded and does not intersect each other. This derives a
contradiction and Lemma \ref{lem2.1} follows. \qed

\begin{la}\label{lem2.1+}\ Let $g(z)$ be an univalent analytic mapping from
$B(a,\eta)$ onto $V$ and $f(w)$ be the inverse of $g(z)$ from $V$
onto $B(a,\eta)$. Then (1) if $0\not\in V$, we have
$$f^\times(g(a))\geq \frac{\eta}{2(1+|a|^2)};$$
(2) if $0\in V$, we have for $|g(a)|\geq L$,
$$f^\times(g(a))\geq \frac{\eta}{4(1+|a|^2)}$$
and for $|g(a)|<L$,
$$f^\times(g(a))\geq \frac{\eta}{2(1+|a|^2)}(\sqrt{1+L^2}-L),$$
where $L=\inf\{|w|: w\not\in V\}.$\end{la}

{\bf Proof.}\ In view of the Koebe quarter covering theorem, we have
$\frac{1}{4}|g'(a)|\eta\leq|g(a)|$ for $0\not\in V$ and
$\frac{1}{4}|g'(a)|\eta\leq|g(a)|+L$ for $0\in V$. Then if $0\not\in
V$, we have
$$f^\times(g(a))=\frac{1+|g(a)|^2}{
|g'(a)|(1+|a|^2)}\geq\frac{\eta}{4(1+|a|^2)}\frac{1+|g(a)|^2}{
|g(a)|} \geq \frac{\eta}{2(1+|a|^2)}.$$ If $0\in V$ and $|g(a)|\geq
L$, we have
$$f^\times(g(a))\geq\frac{\eta}{4(1+|a|^2)}\frac{1+|g(a)|^2}{L+|g(a)|} \geq
\frac{\eta}{4(1+|a|^2)}.$$ If $0\in V$ and $|g(a)|<L$, noticing that
the function $(1+x^2)(L+x)^{-1}$ assumes the minimum value for $x>0$
at $x_L=-L+\sqrt{1+L^2}$, it follows that
$$f^\times(g(a))\geq\frac{\eta}{4(1+|a|^2)}\frac{1+x_L^2}{L+x_L} \geq
\frac{\eta}{2(1+|a|^2)}(\sqrt{1+L^2}-L).$$ \qed

From the proof of Theorem \ref{thm2.2}, in terms of Lemma
\ref{lem2.1} and Lemma \ref{lem2.1+} we have the following

\begin{thm}\label{thm2.3}\ \ Let $f(z)$ be hyperbolic on the Riemann sphere, that is,
$f\in\mathcal{H}(\widehat{\mathbb{C}})$. Then there exist $c>0$ and
$\lambda>1$ such that
\begin{equation}\label{2.2}
(f^n)^\times(z)>c\lambda^n\end{equation} for all $n\geq 1$ and for
all point $z\in \mathcal{J}_f\setminus\mathcal{J}_f(\infty)$. And
for an absolute constant $K$ and for each $n\in\mathbb{N}$, we have
\begin{equation}\label{2.3}d_\infty(f^n(z),f^n(w))\geq K^{-1}c\lambda^n d_\infty(z,w)\end{equation}
whenever $z$ and $w$ are in a common component of
$f^{-n}(B_\infty(a,\eta))$ for some $a\in \widehat{\mathcal{J}}_f$.
\end{thm}

{\bf Proof.}\ \ Set $\alpha(U)=\inf\{f^\times(z):\ z\in U\}$ for an
open set $U$ containing $\mathcal{J}_f$. Since
$f\in\mathcal{H}(\widehat{\mathbb{C}})$, $\alpha>0$ for a $U$. Below
we seek this $U$. Choose a $R>0$ such that $\mathcal{P}(f)\subset
B(0,R)$ and for a sufficiently large $r>R$, all but at most finitely
many components of $f^{-1}(\mathcal{D}_r)$ lie in
$\mathcal{D}_{4R}$. Hence we can assume that each component of
$f^{-1}(\mathcal{D}_r)$ does not intersect $\mathcal{P}(f)$. For a
component $W$ of $f^{-1}(\mathcal{D}_r)$ lying in
$\mathcal{D}_{4r}$, we have at least a point $z_0$ in $W$ at which
$f^\times(z_0)\geq 1$ (Actually, there exists a repelling
fixed-point of $f(z)$ in $W$) and in view of the Koebe distortion
theorem, we have $f^\times(z)\geq K$ for all $z\in W$ where $K$ is
an absolute constant. In terms of Lemma \ref{lem2.1},
$f^{-1}(\mathcal{D}_r)$ has only finitely many components which are
not in $\mathcal{D}_{4r}$ and then there exists a positive constant
$c$ such that $f^\times(z)\geq c$ for all $z\in
f^{-1}(\mathcal{D}_r)$. Now choose finitely many
$a_j\in\mathcal{J}_f\ (1\leq j\leq q)$ such that $B(a_j,\eta)\
(1\leq j\leq q)$ covers $\mathcal{J}_f\cap B(0,2r)$. Since
$f^{-1}(B(a_j,\eta))$ $(j=1,2,...,q)$ has at most one component
which contains $0$ in it, it is obvious in view of Lemma
\ref{lem2.1+} that $f^\times(z)>c_1$ for all $z\in
\cup_{j=1}^qf^{-1}(B(a_j,\eta))$ and for a positive constant $c_1$.
Put $U=\cup_{j=1}^qf^{-1}(B(a_j,\eta))\cup f^{-1}(\mathcal{D}_r)$
and then $\mathcal{J}_f\subset U$ and $\alpha=\alpha(U)>0$.

If $\alpha>1$, then it is obvious that (\ref{2.2}) holds. Assume
that $\alpha\leq 1$. As in the proof of the implication (1)
$\Rightarrow$ (3) of Theorem \ref{thm2.2}, we have a positive
integer $N$ such that $(f^N)^\times(z)>2$ for $z$ inside an open set
$V$ containing $\mathcal{J}_f\setminus \mathcal{J}_f(\infty)$ and we
can assume $V\subseteq U$. For arbitrary positive integer $n$, we
can write $n=mN+r$ for an integer $m$ and $0\leq r<N-1$. Thus for
$z\in \mathcal{J}_f\setminus \mathcal{J}_f(\infty)$, we have
$$(f^n)^\times(z)\geq 2^m(f^r)^\times(z)\geq\alpha^r2^{-1}(2^{1/N})^n\geq c\lambda^n$$
where $\lambda=2^{1/N}$ and $c=2^{-1}\alpha^{N-1}$. This yields
(\ref{2.2}).

To prove the inequality (\ref{2.3}). For each pair of $z, w$ in a
component of $f^{-n}(B_\infty(a,\eta))$ with
$a\in\widehat{\mathcal{J}}_f$, we take a point
$b\in\mathcal{J}_f\setminus \mathcal{J}_f(\infty)$ such that $b\in
B_\infty(a,\eta/2), B_\infty(a,\eta)\subset
B_\infty(b,3\eta/2)\subset B_\infty(b,4\eta)$ and then, in view of
the Koebe distortion theorem, we obtain
\begin{eqnarray*}d_\infty(z,w)&=&d_\infty(f^{-n}_z(f^n(z)),f^{-n}_z(f^n(w)))\\
&\leq& K(f^{-n}_z)^\times(b)d_\infty(f^n(z),f^{n}(w))\\
&=&K\frac{1}{(f^n)^\times(f^{-n}_z(b))}d_\infty(f^n(z),f^{n}(w))\\
&\leq&\frac{K}{c}\lambda^{-n}d_\infty(f^n(z),f^{n}(w))\end{eqnarray*}
where $f^{-n}_z$ is the branch of $f^{-n}$ on $B_\infty(a,\eta)$
which sends $f^n(z)$ to $z$, and equivalently we have (\ref{2.3}).

Thus we complete the proof of Theorem \ref{thm2.3}.\qed

For a hyperbolic meromorphic function with respect to the Euclidean
metric, by the similiar arguments to that in the proof of Theorem
\ref{thm2.3} and in the implication $(3)\Rightarrow (1)$ in the
proof of Theorem \ref{thm2.2}, we can establish the following, which
confirms the expanding property of such meromorphic functions with
respect to the Euclidean metric.

\begin{thm}\label{thm2.3+}\ A meromorphic function $f(z)$ is hyperbolic with respect to the Euclidean
metric, that is, in $\mathcal{H}(Eu)$, if and only if there exist a
$c>0$, a $\lambda>1$ and $\delta>0$ such that for each $n\in
\mathbb{N}$, we have
\begin{equation}\label{2.3+}d(f^n(z),f^n(w))\geq c\lambda^n d(z,w)\end{equation} whenever
$z$ and $w$ are in a common component of $f^{-n}(B(a,\delta))$ for
some $a\in \mathcal{J}_f$.
\end{thm}

It is obvious that (\ref{2.3+}) implies (\ref{1.3}) from the
definition of derivatives. We remark that the $\delta$ in Theorem
\ref{thm2.3+} can be taken to be
$\frac{1}{4}d(\mathcal{P}(f),\mathcal{J}_f)$ when $f$ is proved to
be in $\mathcal{H}(Eu)$.

However, we do not know if the only inequality (\ref{2.2}) implies
the hyperbolicity of the function in question on the Riemann sphere.
It is well known that this thing is true for rational case
(Actually, (\ref{2.2}) implies that the function has no critical
values and no indifferent periodic points on its Julia set).

A hyperbolic function on the complex plane also has some expanding
property on its Julia set, which was proved Rippon and Stallard
\cite{RipponStallard}, that for each $n$ and for each analytic point
$z\in\mathcal{J}_f$ of $f^n$, we have
\begin{equation}\label{2.5+}|(f^n)'(z)|>c\lambda^n\frac{1+|f^n(z)|}{
1+|z|},\end{equation} where $c>0$ and $\lambda>1$ are constants.

In Theorem \ref{thm2.2}, we have known that a hyperbolic function in $\mathcal{H}(\widehat{\mathbb{C}})$
has only finitely many attracting Fatou components and their
preimages under iterates, but no others
(Please notice that this result is also true for hyperbolic function on the complex plane), and therefore the dynamical
behaviors of such functions are clear on their Fatou sets. The
remainder of this section is devoted to study of the Julia set of a
hyperbolic function.

We shall use symbolic dynamics to describe Julia sets when they are
a Cantor set as did in \cite{DevaneyKeen}. Set
$$\Sigma=\{1,2,...\}^\mathbb{N}\cup\{(s_1,s_2,...,s_{n},\infty):
n\in\mathbb{N}, s_j\in\mathbb{N} (j=1,2,...,n)\}.$$ We consider the
topology on $\Sigma$ which is defined in \cite{Moser} as follows: if
$s=(s_1,s_2,...)\in\{1,2,...\}^\mathbb{N}$, then the sets
$$V_k=\{(s_1,s_2,...,s_k,t_{k+1},...): t_j\in\mathbb{N}\ {\rm
for}\ j\geq k+1\}$$ are a neighborhood basis of $s$; if
$s=(s_1,s_2,...,s_{n},\infty)$, then the sets $V_k$ for $k<n$ and
$$W_k=\{(s_1,s_2,...,s_{n},t_{n+1},...): t_{n+1}\geq k\}$$
are a neighborhood basis of $s$. The shift automorphism
$\sigma:\Sigma\rightarrow\Sigma$ is defined by the formula
$\sigma((s_1,s_2,...))=(s_2,...)$ and $\sigma$ is continuous in the
above Moser's topology.

\begin{thm}\label{thm2.4}\ \ Let $f(z)$ be in $\mathcal{H}(\widehat{\mathbb{C}})$ and
the derived set of $\widehat{\rm sing}(f^{-1})$ is finite. If
$\widehat{\mathcal{P}}(f)$ is contained in a component of
$\mathcal{F}_f$, or $\mathcal{F}_f$ is connected, then
$\mathcal{J}_f$ is a Cantor set and $f|_{\mathcal{J}_f}$ is
topologically conjugate to the shift automorphism
$\sigma|_\Sigma$.\end{thm}

{\bf Proof.}\ \ Assume that $\widehat{\mathcal{P}}(f)$ is in a
component $W$ of $\mathcal{F}_f$. Then $W$ is an immediate
attracting domain of $f(z)$ and assume that $a$ is the fixed-point
of $f(z)$ in $W$. Draw a disk $B(a,r)$ in $W$. Since
$\#(\widehat{\rm sing}(f^{-1}))'<\infty$, it is easily seen that
$\#((\widehat{\mathcal{P}}(f))'\setminus B(a,r))<\infty$ and thus we
can draw finitely many disjoint disks $B(a_j,r_j)\ (1\leq j\leq q)$
in $W$ which cover $\widehat{\mathcal{P}}(f)\setminus B(a,r)$.
Adding $q$ disjoint curves $\gamma_j$ in $W$ connecting $B(a_j,r_j)$
and $B(a,r)$, we have
$$U=\widehat{\mathbb{C}}\setminus\left(\bigcup_{j=1}^q(\overline{B}(a_j,r_j)\cup\gamma_j)
\cup \overline{B}(a,r)\right)$$ is a simply connected neighborhood
of $\widehat{\mathcal{J}}_f$ on $\widehat{\mathbb{C}}$ and
$U\cap\widehat{\mathcal{P}}(f)=\emptyset$. It is clear that
$\mathcal{J}_f\subset f^{-1}(U)$ and $f^{-1}(U)$ consists of
infinitely many components all of which have boundaries in the Fatou
set. Thus $\widehat{\mathcal{J}}_f$ is disconnected and $\infty$ is
a single-pointed component of $\widehat{\mathcal{J}}_f$, as
$\infty\not\in f^{-1}(U)$, and so is every point of
$\mathcal{J}_f(\infty)$. In view of the expanding property of $f(z)$
on $\mathcal{J}_f\setminus \mathcal{J}_f(\infty)$, for some $N$,
$f^{-N}(U)$ has a bounded component $V$. Obviously,
$V\cap\mathcal{J}_f\setminus \mathcal{J}_f(\infty)\not=\emptyset$
and $\partial V\subset\mathcal{F}_f$. Furthermore, the diameters of
components of $f^{-n}(U)\cap V$ tend to zero as
$n\rightarrow\infty$. This together with the fact that
$V\cap\mathcal{J}_f\setminus \mathcal{J}_f(\infty)\subset
f^{-n}(U)\cap V$ and $\partial f^{-n}(U)\subset \mathcal{F}_f$
yields that every component of $V\cap\mathcal{J}_f$ is a single
point. By noting that for some $M$,
$\widehat{\mathcal{J}}_f=f^M(V\cap\mathcal{J}_f)$, we have proved
that $\widehat{\mathcal{J}}_f$ is totally disconnected.

Let $\{f_k\}$ be the sequence of all analytic branches of $f^{-1}$
on $U$ and set $D_k=f_k(U)$. $D_k$ is simply connected on
$\widehat{\mathbb{C}}$ and for a pair of distinct $n$ and $m$,
$D_n\cap D_m=\emptyset$ and $\mathcal{J}_f\subset \cup_{k=1}^\infty
D_k$. Define a mapping $\phi:\mathcal{J}_f\rightarrow\Sigma$ as
follows: for a point $z_0\in \mathcal{J}_f$, we determine a point
$s=(s_1,s_2,...,s_n,\ldots)$ on $\Sigma$ by letting $s_n=k$ if
$f^{n-1}(z_0)\in D_k$; $s_n=\infty$ if $f^{n-1}(z_0)=\infty$ and in
this case we shall stop our step at the $n$th entry of $s$.
Obviously, $s$ is uniquely determined by $z_0$. Then set
$s=\phi(z_0)$.

We claim that $\phi(z_0)\not=\phi(z_1)$ for $z_0\not=z_1$. Suppose
that there exist two distinct points $z_0$ and $z_1$ such that
$$\phi(z_0)=\phi(z_1)=(s_1,s_2,...)=s\ ({\rm say}).$$
If $s\in\{1,2,...\}^\mathbb{N}$, then for some $m$,
$f^{-m}_{z_0}(D_{s_{m+1}})\cap f^{-m}_{z_1}(D_{s_{m+1}})=\emptyset$
where $f^{-m}_{z_i}\ (i=0,1)$ are respectively the branches of
$f^{-m}$ with $f^{-m}_{z_i}(f^m(z_i))=z_i.$ Notice that $f^m(z_i)\in
D_{s_{m+1}}\cap U\ (i=0,1)$ and
$D_{s_{m+1}}\cap\mathcal{P}(f)=\emptyset$ and hence every branch of
$f^{-m}$ can be analytically extended to the domain $D_{s_{m+1}}\cup
U$. Consider the branches $f^{-1}_{0,m}$ and $f^{-1}_{1,m}$ of
$f^{-1}$ on $U$ such that $f^{-1}_{0,m}(f^m(z_0))=f^{m-1}(z_0)$ and
$f^{-1}_{1,m}(f^m(z_1))=f^{m-1}(z_1)$. Since both of $f^{m-1}(z_i)\
(i=0,1)$ are in $D_{s_m}$,
$f^{-1}_{0,m}(z)=f_{s_m}(z)=f^{-1}_{1,m}(z)$. Then inductively, we
shall have
$$f^{-m}_{z_0}=f^{-1}_{0,1}\circ...\circ f^{-1}_{0,m}
=f^{-1}_{1,1}\circ...\circ f^{-1}_{1,m}=f^{-m}_{z_1}.$$ This derives
a contradiction. If $s=(s_1,s_2,...,s_n,\infty)$, then we have
$$z_0=f^{-1}_{s_1}\circ...\circ f^{-1}_{s_n}(\infty)=z_1,$$ but it
contracts the assumption of $z_0\not=z_1$. Thus we complete the
proof of our claim.

For each $k$, define a set mapping $\widetilde{f}_k$: for every
$D\subseteq U$, $\widetilde{f}_k(D)=f_k(D)\cap U$. For a
$s=\{s_1,s_2,...,s_n,\ldots\}$ on $\Sigma$ such that each
$s_n\not=\infty$, according to the expanding property, the diameters
of $\widetilde{f}_{s_1}\circ \widetilde{f}_{s_2}\circ\ldots\circ
\widetilde{f}_{s_n}(U)$ tend to zero as $n\rightarrow\infty$ and
since $\widetilde{f}_{s_n}(U)\subset U$, we have
$\widetilde{f}_{s_1}\circ \widetilde{f}_{s_2}\circ\ldots\circ
\widetilde{f}_{s_n}(U)\subset \widetilde{f}_{s_1}\circ
\widetilde{f}_{s_2}\circ\ldots\circ \widetilde{f}_{s_{n-1}}(U)$ and
each $\widetilde{f}_{s_1}\circ \widetilde{f}_{s_2}\circ\ldots\circ
\widetilde{f}_{s_n}(U)$ intersects
$\mathcal{J}_f\setminus\mathcal{J}_f(\infty)$. Therefore, the set
$\bigcap_{n=1}^\infty \widetilde{f}_{s_1}\circ
\widetilde{f}_{s_2}\circ\ldots\circ \widetilde{f}_{s_{n}}(U)$ is a
single point $z_0\in\mathcal{J}_f\setminus\mathcal{J}_f(\infty)$. We
claim that $s=\phi(z_0)$. For each $n\geq 1$, $z_0\in
\widetilde{f}_{s_1}\circ \widetilde{f}_{s_2}\circ\ldots\circ
\widetilde{f}_{s_{n}}(U)=
f_{s_1}(\widetilde{f}_{s_2}\circ\ldots\circ
\widetilde{f}_{s_{n}}(U))\cap U$, that is,
$z_0\in\widetilde{f}_{s_2}\circ\ldots\circ
\widetilde{f}_{s_{n}}(U)$. In general, we have $f^{n-1}(z_0)\in
\widetilde{f}_{s_{n}}(U)= f_{s_{n}}(U)\cap U.$ This implies that
$f^{n-1}(z_0)\in D_{s_n}$. According to the definition of $\phi$, we
have proved the claim. For a $s=\{s_1,s_2,...,s_n,\infty\}$ on
$\Sigma$, $f_{s_1}\circ f_{s_2}\circ\ldots\circ f_{s_{n}}(\infty)$
is a single point $z_0\in\mathcal{J}_f(\infty)$ and $s=\phi(z_0)$.

It is easy to prove that $\phi$ is a homeomorphism from
$\mathcal{J}_f$ onto $\Sigma$ and $\phi\circ
f(z)=\sigma\circ\phi(z)$ on $\mathcal{J}_f$. Therefore
$\widehat{\mathcal{J}}_f$ is a Cantor set and we have completed
proof of Theorem \ref{thm2.4}.\qed

Theorem \ref{thm2.4} was proved by Steinmetz in his book
\cite{Steinmetz} for hyperbolic rational functions. For a
meromorphic function $f(z)$, if its post-singular set is contained
in a Fatou component $W$ of it, then $W$ is completely invariant and
an attracting domain of $f(z)$. According to the relation of Fatou
components with post-singular set, it is easily seen that $f(z)$ has
no other Fatou components than $W$ and so $\mathcal{F}_f=W$, that
is, $\mathcal{F}_f$ is connected.

The Julia set of a hyperbolic function in
$\mathcal{H}(\widehat{\mathbb{C}})$ may not be totally disconnected
even if it is disconnected, which is explained by the function
$R(z)=z^2+\lambda/z^3$, but no such transcendental example has been
found. Indeed, it was proved by McMullen \cite{McMullen} that the
Julia set of $R(z)$ is a Cantor set of circles for sufficiently
small $\lambda>0$ and $\infty$ attracts all critical points of
$R(z)$ and so $R(z)$ is hyperbolic in
$\mathcal{H}(\widehat{\mathbb{C}})$. If Julia set of a meromorphic
function $f(z)$ which is not of the form
$\alpha+(z-\alpha)^{-k}e^{g(z)}$ for a natural number $k$, a complex
number $\alpha$ and an entire function $g(z)$, is disconnected on
$\widehat{\mathbb{C}}$, then it has uncountably infinitely many
Julia components and it was proved in Ng, Zheng and
Choi\cite{NgZheng} that it has uncountably infinitely many buried
components if $\mathcal{F}_f$ has no completely invariant
components. Since $\mathcal{J}_f(\infty)$ has at most countably many
points, $\mathcal{J}_f$ has only countably components which contain
points of $\mathcal{J}_f(\infty)$ and if it is disconnected on
$\widehat{\mathbb{C}}$, then $\widehat{\mathcal{J}}_f$ has
uncountably infinitely many components which do not contain any
points in $\mathcal{J}_f(\infty)$. It is clear that $f(z)$ maps a
component of $\mathcal{J}_f\setminus\mathcal{J}_f(\infty)$ into a
component of $\mathcal{J}_f\setminus\mathcal{J}_f(\infty)$ (if
$f(z)$ is rational, we use $\mathcal{J}_f$ in the place of
$\mathcal{J}_f\setminus\mathcal{J}_f(\infty)$). A component
$\mathcal{J}_0$ of $\mathcal{J}_f\setminus\mathcal{J}_f(\infty)$ of
$f(z)$ is called periodic if for some $n\geq 1$,
$f^n(\mathcal{J}_0)\subseteq \mathcal{J}_0$ and the smallest $n$
with this property is period of $\mathcal{J}_0$; pre-periodic if for
some $n\geq 1$, the component containing $f^n(\mathcal{J}_0)$ is
periodic, but $\mathcal{J}_0$ is not periodic; wandering if for
$n\not=m$, $f^n(\mathcal{J}_0)\cap f^m(\mathcal{J}_0)=\emptyset.$ A
component of $\mathcal{J}_f\setminus\mathcal{J}_f(\infty)$ may not
be closed on $\widehat{\mathbb{C}}$. The transcendental case is more
complicated than the rational case and the complicated mainly
results from the essential singular point $\infty$. Let us explain
that. For a continuum component $\mathcal{J}_1$ of $\mathcal{J}_f$
containing a pole of $f(z)$, $f(\mathcal{J}_1)$ is connected on
$\widehat{\mathbb{C}}$, but may not be connected on $\mathbb{C}$;
If, in addition, there exist two finite asymptotic values $a$ and
$b$ of $f(z)$ on two distinct components of $\mathcal{J}_f$, then
there exist at least two components $\mathcal{J}_1^1$ and
$\mathcal{J}_1^2$ of $f(\mathcal{J}_1)$ on $\mathbb{C}$ such that
$f(\mathcal{J}_1^1)$ and $f(\mathcal{J}_1^2)$ lie respectively in
the components containing $a$ and $b$ and therefore
$f(\mathcal{J}_1^1)$ and $f(\mathcal{J}_1^2)$ do not intersect each
other on $\widehat{\mathbb{C}}$, that is, we mapped a component of
$\mathcal{J}_f$ into several disjoint components of $\mathcal{J}_f$.
However, the situation cannot happen if we consider components of
$\mathcal{J}_f\setminus\mathcal{J}_f(\infty)$, but we pay a price
for that a component of $\mathcal{J}_f$ could be divided into
several components of $\mathcal{J}_f\setminus\mathcal{J}_f(\infty)$.
Here is an example. If $\mathcal{J}_f$ contains a locally isolate
Jordan arc, then a result of Stallard \cite{Stallard3} yields that
$\mathcal{J}_f$ is a Jordan curve or arc and in this case, since
$\mathcal{J}_f(\infty)$ is dense in $\mathcal{J}_f$, every component
of $\mathcal{J}_f\setminus\mathcal{J}_f(\infty)$ is single-pointed.

McMullen in \cite{McMullen} proved that the Julia set of a rational
function has at most countably many periodic or preperiodic components. We do
not know if the result is true for the transcendental case. If the
Julia set is totally disconnected, then the question is affirmative.
Pilgrim and Tan characterized in \cite{PilgrimTan} Julia components
of a hyperbolic rational function by proving that with the possible
exception of finitely many periodic components and their countable
collection of preimages, every connected component of the Julia set
of a hyperbolic rational function is either a point or a Julia
curve. This leads us to pose the following question.

\begin{que}\label{que2.1} Is every connected component of the Julia set of a transcendental
hyperbolic meromorphic function in $\mathcal{H}(\widehat{\mathbb{C}})$ either a point or a Julia curve with
possible exception of at most countably infinitely many periodic or
pre-periodic components or those components which contain points of
$\mathcal{J}_f(\infty)$?\end{que}

If a Fatou component of a hyperbolic function is multiply-connected
with connectivity at least three, then the single-pointed components of its Julia set are
dense in its Julia set (cf. Dominguez \cite{Dominguez}). However, in
this case we do not know if the Julia set is totally disconnected,
which produces the following.

\begin{que}\label{que2.2} Is the Julia set of a hyperbolic meromorphic
function in $\mathcal{H}(\widehat{\mathbb{C}})$ a Cantor set if it has a single-pointed component or the
function has a multiply-connected Fatou component with connectivity
at least three?
\end{que}

If the Julia set of a hyperbolic function is connected on
$\widehat{\mathbb{C}}$, then what can we say? It is well known that
a hyperbolic rational function has a locally connected Julia set if
the Julia set is connected. By the same method as in the proof of
rational case, we can establish the result that if $f(z)$ is a
hyperbolic transcendental meromorphic function in
$\mathcal{H}(\mathbb{C})$ and $\mathcal{J}_f$ is connected and if
$\mathcal{F}_f$ has only bounded components, that is, $\infty$ is a
buried point of $\mathcal{J}_f$, then $\mathcal{J}_f$ is locally
connected. In fact, since $\mathcal{F}_f$ has only bounded
components, no poles lie on the boundaries of Fatou components and
$f(z)$ has no asymptotic values (Indeed, if $f(z)$ has an asymptotic
value, then the asymptotic value must be in the Faout set and hence
there exists at least a tract corresponding to the value which is
contained in the Fatou set. This implies that the Fatou set has an
unbounded component). Therefore $f(z)$ is a proper mapping from a
Fatou component onto a Fatou component and the dynamical behavior of
$f(z)$ on an attracting domain is conjugate to a finite Blaschke
product on the unit disk. From this it follows that the boundary of
the attracting domain is a Jordan curve. For hyperbolic case in
$\mathcal{H}(\widehat{\mathbb{C}})$ we pose the following

\begin{que}\label{que2.3} Is the Julia set of a hyperbolic meromorphic
function in $\mathcal{H}(\widehat{\mathbb{C}})$ locally connected if it is connected on $\widehat{\mathbb{C}}$? \end{que}

Generally, Question \ref{que2.3} is negative for a function in
$\mathcal{H}(\mathbb{C})$. This is because an unbounded attracting
domain of a transcendental entire function has boundary which is not
locally connected (see Baker and Domingueze \cite{BakerDomingueze}).
The function $\lambda\sin z$ with $|\lambda|<1$ is in
$\mathcal{H}(\mathbb{C})$ and its Julia set is connected on
$\widehat{\mathbb{C}}$, but not locally connected at $\infty$.
Actually, its Fatou set is an unbounded attracting domain and
$\infty$ is inaccessible in the Fatou set. In view of the same
argument as in Bergweiler and Eremenko \cite{BergweilerEremenko}, we
can prove the following

\begin{thm}\label{thm2.5}\ Let $f(z)$ be a hyperbolic meromorphic
function in $\mathcal{H}(\mathbb{C})$. If the Fatou set $f(z)$
consists exactly of two completely invariant components, then the
Julia set $\mathcal{J}_f$ is a Jordan curve.\end{thm}

It is a natural thinking that if the Julia set of a meromorphic
function is simple, then the behavior of it near $\infty$ should be
simple. We have the following

\begin{thm}\label{thm2.6}\ \ Let $f(z)$ be a meromorphic
function in Class $\mathcal{B}$, that is, ${\rm sing}(f^{-1})$ is bounded, with connected and locally connected
Julia set. If the Fatou set of $f(z)$ has exactly two components,
then $\infty$ is not an asymptotic value of $f(z)$.\end{thm}

{\bf Proof.}\ \ Take a $R>0$ such that ${\rm sing}(f^{-1})\subset
\mathbb{C}\setminus B(0,R)$ and draw a Jordan curve $\gamma$ outside
$B(0,R)$ which surrounds the origin and intersects $\mathcal{J}_f$ at only two points $a$ and
$b$. This can be done because $\mathcal{J}_f$ is connected and
locally connected. If $\infty$ is an asymptotic value of $f(z)$,
then the inverse of $f(z)$ has a logarithmic singularity over
$\infty$. Let $U$ be the component of $f^{-1}({\rm out}(\gamma))$
such that $f:U\rightarrow {\rm out}(\gamma)$ is a universal
covering. Since $\mathcal{J}_f\cap{\rm out}(\gamma)$ has two
components on $\mathbb{C}$, there exist at least two unbounded
components $\gamma_0$ and $\gamma_1$ of $U\cap \mathcal{J}_f$ such
that $\gamma_0$ starts from a point $z_0$ of $f^{-1}(a)$ and
$\gamma_1$ from a point $z_1$ of $f^{-1}(b)$. Since the Fatou set of
$f(z)$ consists of two components $V$ and $W$, these two components
$V$ and $W$ are completely invariant under $f(z)$ or $f^2(z)$ and
hence $f(z)$ must have infinitely many poles and
$\mathcal{J}_f=\partial V=\partial W$. It is clear that each pole of
$f(z)$ is outside $U$ ($\infty$ is a logarithmic singular value of $f(z)$). And since $\partial U$ does not wind around
$\infty$, $\mathcal{J}_f\cap(\mathbb{C}\setminus U)$ has an unbounded
component $\gamma$ connecting $z_0$ and $z_1$ and thus
$\mathcal{F}_f$ has at least three components. This contradicts our
assumption, although these three components may have the common
boundary. We have proved that $\infty$ is not an asymptotic value of
$f(z)$.\qed

In view of Theorem \ref{thm2.6}, the function satisfying the
assumption of Theorem \ref{thm2.5} is actually in
$\mathcal{H}(\widehat{\mathbb{C}})$ by noticing that $\infty$ cannot
be a critical value as the Julia set is a Jordan curve. Combining
Theorem of \cite{BergweilerEremenko} and Theorem \ref{thm2.6} yields
the following

\begin{cor}\label{cor2.1}\ Let $f(z)$ be a meromorphic function in Class $\mathcal{S}$, that is,
${\rm sing}(f^{-1})$ is finite. If $f(z)$ has two completely invariant components, then
$\infty$ is not a singular value of $f(z)$. Furthermore, if, in addition, $f(z)$ has no a fixed point with multiplier
$1$, then $f(z)$ is hyperbolic in $\mathcal{H}(\widehat{\mathbb{C}})$. \end{cor}

{\bf Proof.}\ In view of Theorem of \cite{BergweilerEremenko}, $\widehat{\mathcal{J}}_f$ is a Jordan curve and therefore,
$\infty$ is not a critical value of $f(z)$. By means of Theorem \ref{thm2.6}, $\infty$ is not an asymptotic value of $f(z)$.
Thus $\infty$ is not a singular value of $f(z)$. Using Theorem of \cite{BergweilerEremenko}, we have
$\widehat{{\rm sing}}(f^{-1})={\rm sing}(f^{-1})\subset \mathcal{F}_f$ and $\widehat{\mathcal{P}}(f)\subset \mathcal{F}_f$ if
$f(z)$ has no a fixed point with multiplier
$1$. Thus $f(z)$ is in $\mathcal{H}(\widehat{\mathbb{C}})$. \qed

\vskip 1cm
\section{Dimension, Conformal and Invariant Measures for Hyperbolic Functions}

This section is devoted to the discussion of existence of conformal
invariant measure on the Julia set and the Hausdorff dimension of
the Julia set of a hyperbolic function. The results we shall obtain
are well-known for a rational hyperbolic function. First of all, we
give out a result about area of the Julia set.

\begin{thm}\label{thm3.1}\ \ Let $f(z)$ be a transcendental
meromorphic function such that
\begin{equation}\label{3.1}\#(\mathcal{J}_f\cap\mathcal{P}(f))<\infty\end{equation}
and $\infty$ is not a singular value of $f(z)$. Assume that
$\mathcal{P}(f)\cap\mathcal{J}_f(\infty)=\emptyset$. Then
${\rm Area}(\mathcal{J}_f)=0.$ In particular, if $f(z)$ is
hyperbolic in $\mathcal{H}(\widehat{\mathbb{C}})$, then ${\rm Area}(\mathcal{J}_f)=0.$
\end{thm}

The function $f(z)$ satisfying (\ref{3.1}) is called geometrically
finite on the complex plane and if, in addition, $\infty$ is not an
asymptotic value of $f(z)$, then $f(z)$ is called geometrically
finite on the Riemann sphere. Obviously a hyperbolic function on the
Riemann sphere (resp. on the complex plane) is geometrically finite
on the Riemann sphere (resp. on the complex plane) and hence in view
of the former half part of Theorem \ref{thm3.1}, the Julia set of a
hyperbolic function on the Riemann sphere has zero area. However,
the result is not true for a hyperbolic function on the complex
plane. The condition
"$\mathcal{P}(f)\cap\mathcal{J}_f(\infty)=\emptyset$" implies that
$\infty$ is not in $\mathcal{P}(f)$ and so $\mathcal{P}(f)$ is
compact on $\mathbb{C}$. From this together with the assumption that
$\infty$ is not a singular value of $f(z)$, it follows that
$\mathcal{J}_f$ is thin at $\infty$ and then Theorem \ref{thm3.1}
can be proved by a result of Zheng \cite{Zheng3}.

Sullivan \cite{Sullivan} investigated in terms of the derivative
with respect to the sphere metric the Hausdorff dimension of the
Julia set of a hyperbolic rational function $f(z)$ and conjectured
that the Hausdorff dimension depends real analytically on $f$, which
was proved by Ruelle \cite{Ruelle} in view of the Bowen's formula.
Many mathematicians investigated the case of transcendental
meromorphic functions, please see Baranski \cite{Baranski}, Kotus
and Urbanski \cite{KotusUrbanski}, Mayer and Urbanski
\cite{MayerUrbanski} and Stallard \cite{Stallard1} and so on.

In what follows, we discuss the thermodynamic formalism of
hyperbolic transcendental meromorphic functions. Let $f(z)$ be a
transcendental meromorphic function with
$\widehat{\mathcal{J}}_f\cap\mathcal{P}(f)=\emptyset$, that is,
$f(z)$ is hyperbolic on the complex plane. For $\varphi\in
C(\mathcal{J}_f\setminus f^{-1}(\infty))$, define the
Perron-Frobenius-Ruelle operator $\mathcal{L}_\varphi$ on
$C(\widehat{\mathcal{J}}_f)$ by the formula
$$\mathcal{L}_\varphi(g)(a)=\sum_{f(z)=a}g(z)e^{\varphi(z)}$$ for
$a\in \mathcal{J}_f$. Then a simple calculation yields
$$\mathcal{L}^n_\varphi(g)(a)=\mathcal{L}_\varphi(\mathcal{L}^{n-1}_\varphi(g))(a)=
\sum_{f^n(z)=a}g(z)\exp(S_n\varphi(z))$$ where
$S_n\varphi(z)=\sum_{i=0}^{n-1}\varphi(f^i(z)).$ In particular, for
$t>0$ and $\varphi_t(z)=-t\log f^\times(z)$, we write $\mathcal{L}^n_t$ for $\mathcal{L}^n_{\varphi_t}$ and we have
$$\mathcal{L}^n_t(1)(a)=\mathcal{L}^n_{\varphi_t}(1)(a)=
\sum_{f^n(z)=a}\frac{1}{((f^n)^\times(z))^t}$$ and define
$$P_a(f,t)=\limsup_{n\rightarrow\infty}\frac{1}{n}\log \mathcal{L}^n_t(1)(a).$$
The function $P_a(f,t)$ plays a key role in our discussion. It is
important that for what $t$, $\mathcal{L}^n_t(1)(a)$ or $P_a(f,t)$
is finite. This is the first question we should answer.

Stallard \cite{Stallard} discussed the Hausdorff dimension of a
hyperbolic function on the complex plane and established the
following

\begin{thm}\label{thm3.2}\ Let $f(z)$ be a transcendental
meromorphic function in $\mathcal{H}(\mathbb{C}).$ Then there exists
a real number $s(f)$ with $0<s(f)\leq 2$ such that

(1) for every $a\in\mathcal{J}_f$, we have
\begin{equation}\label{3.2}\lim\limits_{n\rightarrow\infty}
\sum_{f^n(z)=a}\frac{1}{((f^n)^\times(z))^t}=\begin{cases}\infty
& t<s(f),\\ 0, & t>s(f)\end{cases}\end{equation}

(2) ${\rm dim}_H(\mathcal{J}_f)\geq s(f).$\end{thm}

Here and throughout the whole section the notation ${\rm dim}_H(X)$
denotes the Hausdorff dimension of set $X$.

Actually, we have $s(f)=\inf\{t\geq 0: P_a(f,t)\leq 0\}$, which is
called Poincar\'e exponent. From the proof of Lemma 7.3.2 of
\cite{Zheng} we have the following result, whose proof will be given
for completeness.

\begin{la}\label{lem3.3+}\ Let $f(z)$ be a transcendental
meromorphic function in $\mathcal{H}(\mathbb{C}).$ For a point
$a\in\widehat{\mathcal{J}}_f$ and each $n\geq 1$, assume that
$g_n(z)$ is a single-valued analytic branch of $f^{-n}$ on
$B(a,4\delta)$. Then there exists a positive number $\rho$ only
depending on $\delta$ such that for $t>0$, we have a constant
$K_t>0$ such that
\begin{equation}\label{3.3}\left|1-\frac{((f^n)^\times(g_n(u)))^t}{
((f^n)^\times(g_n(v)))^t}\right|\leq K_td(u,v)\end{equation} for
arbitrary two points $u, v\in B(a,\rho)$, where for $a=\infty$, we
use $B_\infty$ and $d_\infty$ in the places of $B$ and $d$ and in
this case, we assume, in addition, that $0\in
\mathcal{F}_f$.\end{la}

{\bf Proof.}\ Since $u, v\in B(a,\delta)$, we have
\begin{eqnarray*}\left|1-\frac{1+|v|^2}{
1+|u|^2}\right|&=&\frac{||v|^2-|u|^2|}{
1+|u|^2}\leq|u-v|\frac{|v|+|u|}{
1+|u|^2}\nonumber\\
&\leq&|u-v|\frac{2|u|+|u-v|}{
1+|u|^2}\leq(1+\delta)d(u,v)\end{eqnarray*} and equivalently,
\begin{equation}\label{3.4+1}\frac{1+|v|^2}{1+|u|^2}=1+C_1d(u,v)\end{equation}
for some real number $C_1$ with $|C_1|\leq 1+\delta.$ In view of the
Koebe distortion theorem, we have
\begin{eqnarray*}\frac{|g_n'(v)|}{
|g_n'(u)|}&\leq&\left(\frac{2\delta+|u-v|}{2\delta-|u-v|}\right)^4\\
&\leq&\exp\left[4\log\left(1+\frac{2}{\delta}|u-v|\right)\right]\\
&\leq&\exp(8\delta^{-1}|u-v|)\\
&\leq&1+e^{16}8\delta^{-1}|u-v|\end{eqnarray*} and
$$\frac{|g_n'(v)|}{|g_n'(u)|}\geq\exp(-8\delta^{-1}|u-v|)\geq 1-8\delta^{-1}|u-v|.$$
Thus it is easily seen that
\begin{equation}\label{3.4+2}\frac{|g_n'(v)|}{
|g_n'(u)|}=1+C_2d(u,v)\end{equation} for some real number $C_2$ with
$|C_2|\leq e^{16}8\delta^{-1}.$

In view of Theorem of \cite{RipponStallard}, for each
$z\in\mathcal{J}_f\setminus\mathcal{J}_f(\infty)$ and some
$\lambda>1$ and $c>0$, we have (\ref{2.5+}). Using the Koebe
distortion theorem together with (\ref{2.5+}) yields
\begin{eqnarray}\label{3.4+3}
|g_n(u)-g_n(v)|&\leq&\max\{|g'_n(z)|:\ \forall z\in
B(a,\delta)\}|u-v|\nonumber\\
&\leq&81|g_n'(v)||u-v|\nonumber\\
&\leq&81|u-v|\frac{1}{|(f^n)'(g_n(v))|}\nonumber\\
&\leq&81|u-v|\frac{1}{c\lambda^n}\frac{|g_n(v)|+1}{
|v|+1}.\end{eqnarray} When $|g_n(v)|<1$, it follows from
(\ref{3.4+3}) that $|g_n(u)-g_n(v)|\leq 162c^{-1}|u-v|$ and so
$|g_n(u)|\leq|g_n(v)|+162c^{-1}|u-v|<1+324c^{-1}\delta$; When
$|g_n(v)|\geq 1$, we have $|g_n(u)-g_n(v)|\leq
162c^{-1}|u-v||g_n(v)|$ and so $|g_n(u)|\leq
(1+162c^{-1}|u-v|)|g_n(v)|\leq (1+324c^{-1}\delta)|g_n(v)|.$ Thus we
always have
\begin{eqnarray*}\left|1-\frac{1+|g_n(u)|^2}{
1+|g_n(v)|^2}\right|&\leq&|g_n(u)-g_n(v)|\frac{|g_n(u)|+|g_n(v)|}{
1+|g_n(v)|^2}\\
&\leq&162c^{-1}|u-v|(1+|g_n(v)|)\frac{(1+324c^{-1}\delta)(1+|g_n(v)|)}{
1+|g_n(v)|^2}\\
&\leq&324c^{-1}(1+324c^{-1}\delta)d(u,v)\end{eqnarray*} and
equivalently, \begin{equation}\label{3.4+4}\frac{1+|g_n(u)|^2}{
1+|g_n(v)|^2}=1+C_3d(u,v)\end{equation} for some real number $C_3$
with $|C_3|\leq 324c^{-1}(1+324c^{-1}\delta).$

Combining the above equalities (\ref{3.4+1}), (\ref{3.4+2}) and
(\ref{3.4+3}) yields that \begin{eqnarray*}\frac{g_n^\times(v)}{
g_n^\times(u)}&=&(1+C_1d(u,v))(1+C_2d(u,v))(1+C_3d(u,v))\\
&=&1+D_1d(u,v)\end{eqnarray*} for some real number $D_1$ with
$$|D_1|\leq
|C_1|+|C_2|+|C_3|+(|C_1C_2|+|C_1C_2|+|C_2||C_3|)2\delta+|C_1C_2C_3|4\delta^2.$$
The quantity on the right side of the above inequality is not larger
than a positive constant $D$ which depends only on $c$ and $\delta$.
Now choose an $\rho>0$ such that $2D\rho\leq\frac{1}{2}$. For $t\geq
1$, we have
\begin{eqnarray*}\left|1-\left(\frac{g_n^\times(v)}{
g_n^\times(u)}\right)^t\right|&=&|1-(1+D_1d(u,v))^t|\\
&=&t|D_1|d(u,v)\int_0^1(1+D_1d(u,v)x)^{t-1}dx\\
&\leq&t2^{t-1}Dd(u,v)\end{eqnarray*} and for $0<t<1$,
\begin{eqnarray*}\left|1-\left(\frac{g_n^\times(v)}{
g_n^\times(u)}\right)^t\right|&=&t|D_1|d(u,v)\int_0^1\frac{1}{(1+D_1d(u,v)x)^{1-t}}dx\\
&\leq&t2^{1-t}Dd(u,v),\end{eqnarray*} whenever $u,v\in B(a,\rho)$.
Thus from the above two inequalities we obtain (\ref{3.3}) for
$K_t=\max\{t2^{t-1}D,t2^{1-t}D\}$.

We can prove the case for $a=\infty$ by using the same argument as
above in terms of the sphere metric instead of the Euclidean metric.
\qed

To discuss the finiteness of $\mathcal{L}^n_t(1)(a)$, we establish
the following

\begin{la}\label{lem3.1+}\ Let $f(z)$ be a transcendental
meromorphic function in $\mathcal{H}(\mathbb{C}).$

(i) There exists a $\rho>0$ such that for each pair of $a,
b\in\mathcal{J}_f$ with $d(a,b)<\rho$ we have
\begin{equation}\label{3.8}\mathcal{L}^n_t(1)(b)\leq
(1+K_td(a,b))\mathcal{L}^n_t(1)(a);\end{equation}

(ii) If for some $n\geq 1$ and some $a\in \mathcal{J}_f$,
$\mathcal{L}^n_t(1)(a)$ is finite, then for each pair of positive
integers $p$ and $q$ with $p+q=n$ and each $w\in f^{-p}(a)$, we have
\begin{equation}\label{3.9}\mathcal{L}^q_t(1)(w)\leq
\mathcal{L}^n_t(1)(a)(f^p)^\times(w)^t\end{equation} and hence if
for all large $n$, $\mathcal{L}^n_t(1)(a)$ is finite, then for each
$p>0$, $\mathcal{L}^p_t(1)(w)$ is finite on $\mathcal{J}_f$;

(iii) For $t>s(f)$ and each $n\geq 1$, $\mathcal{L}^n_t(1)(w)$ is
finite on $\mathcal{J}_f$;

(iv) If ${\rm sup}\{\mathcal{L}_t(1)(a): a\in\mathcal{J}_f\}\leq M$,
then ${\rm sup}\{\mathcal{L}^n_t(1)(a):\ a\in\mathcal{J}_f\}\leq
M^n$ so that $P_a(f,t)\leq M$ for each $a\in\mathcal{J}_f$.\end{la}

{\bf Proof.}\ We can find a fixed $\rho>0$ such that (\ref{3.3})
holds for each $a\in\mathcal{J}_f$ and from this we easily obtain
(\ref{3.8}). And (\ref{3.9}) is obvious from the expression of
$\mathcal{L}^n_t(1)(a)$ and the equality $(f\circ
g)^\times(z)=f^\times(g(z))g^\times(z)$. Actually, we have
\begin{eqnarray*}\mathcal{L}^{p+q}_t(1)(a)&=&\sum_{f^{p+q}(z)=a}\frac{1}{((f^{p+q})^\times(z))^t}\\
&=&\sum_{f^p(c)=a}\sum_{f^{q}(z)=c}\frac{1}{((f^{p})^\times(c))^t((f^q)^\times(z))^t}\\
&\geq&\sum_{f^{q}(z)=w}\frac{1}{((f^{p})^\times(w))^t((f^n)^\times(z))^t}\\
&=&\frac{1}{((f^{p})^\times(w))^t}\mathcal{L}^q_t(1)(w).\end{eqnarray*}
This is (\ref{3.9}). Now assume that for all large $n$,
$\mathcal{L}^n_t(1)(a)$ is finite. For a fixed $w\in \mathcal{J}_f$,
we can find a large $m$ and a point $c\in \mathcal{J}_f$ with
$d(w,c)<\rho$ such that $f^m(c)=a$, and in view of (\ref{3.9}) we
have $$\mathcal{L}^p_t(1)(c)\leq
\mathcal{L}^{m+p}_t(1)(a)(f^m)^\times(c)^t.$$ This implies that
$$\mathcal{L}^p_t(1)(w)\leq(1+K_t\rho)\mathcal{L}^p_t(1)(c)\leq(1+K_t\rho)
\mathcal{L}^{m+p}_t(1)(a)(f^m)^\times(c)^t.$$ Thus we complete the
proof of the result (ii).

The result (iii) follows from the result (ii) and (\ref{3.2}) and
the result (iv) is easily proved by the following inequality
$$\sum_{f\circ g(z)=a}\frac{1}{(f\circ g)^\times(z)^t}\leq
\sup\limits_{c_n}\left\{\sum_{g(z)=c_n}\frac{1}{
g^\times(z)^t}\right\} \sum_{f(z)=a}\frac{1}{f^\times(z)^t},$$ where
$\{c_n\}$ is the sequence of all roots of $f(z)=a.$ \qed

\begin{la}\label{lem3.1}\ Let $f(z)$ be a transcendental
meromorphic function in $\mathcal{H}(\mathbb{C}).$ Then $P_a(f,t)$
does not depend on $a\in\mathcal{J}_f$; So we write $P(f,t)$ for
$P_a(f,t)$; $P(f,t):(\tau(f), 2]\rightarrow \mathbb{R}$ is a
strictly decreasing convex function in $t$ where
$\tau(f)=\inf\{t\geq 0: P(f,t)<+\infty\}$ and when $s(f)>\tau(f)$,
$P(f,s(f))=0$.\end{la}

{\bf Proof.}\ \ Set
$$\delta=\frac{1}{4}d(\mathcal{J}_f,\mathcal{P}(f))>0.$$
For every $a\in \mathcal{J}_f$ and for every point $z_0$ such that
$f^n(z_0)=a$, we have an analytic branch $g_n$ of $f^{-n}$ on
$B(a,4\delta)$ sending $a$ to $z_0$. In view of the Koebe Theorem,
for every $b\in B(a,\delta)$, we have $K^{-1}g^\times_n(b)\leq
g^\times_n(a)\leq Kg^\times_n(b)$ for an absolute constant $K>1$ and
equivalently
$$K^{-1}(f^n)^\times(g_n(b))\leq
(f^n)^\times(g_n(a))\leq K(f^n)^\times(g_n(b)).$$

For a point $b$ in $\mathcal{J}_f$, we have a point $z_1\in
B(b,\delta)$ such that $f^m(z_1)=a$ for some $m\geq 0$. Then in view
of (\ref{3.9}), we have
$$\mathcal{L}^n_t(1)(b)\leq
K^t\mathcal{L}^n_t(1)(z_1)\leq K^t((f^{m})^\times(z_1))^t
\mathcal{L}^{m+n}_t(1)(a).$$ This implies that $P_a(f,t)\geq
P_b(f,t)$ and so $P_a(f,t)=P_b(f,t)$.

It follows immediately from the H\"older inequality that $P(f,t)$ is
convex in $t$ and it is proved in Theorem 6.3.12 of \cite{Zheng1}
that $P(f,t)$ is strictly decreasing in $t$.

It is obvious that $\tau(f)\leq s(f)$ and it is easily seen that for
$t>s(f)$, $P(f,t)\leq 0$ and for $t<s(f)$, $P(f,t)\geq 0$. Since
$P(f,t)$ is continuous at $s(f)$ when $\tau(f)<s(f)$, we have
$P(f,s(f))=0$. \qed

The case ${\rm dim}_H(\mathcal{J}_f)>s(f)$ in the result (2) of
Theorem \ref{thm3.1} is possible, which was shown by Stallard
\cite{Stallard} by observing an example of entire function in
$\mathcal{H}(\mathbb{C})$. We shall give out another simple example
late for that. Therefore, under what additional condition does ${\rm
dim}_H(\mathcal{J}_f)$ equal to $s(f)$? For this question, we
consider a subset of the Julia set,that is, so-called radial (or
conical) Julia set. For a meromorphic function $f(z)$, a point
$z_0\in\mathcal{J}_f$ is called conical if all forward images
$f^n(z_0)$ of $z_0$ are well defined and there is some
$\delta(z_0)>0$ such that for infinitely many $n\in\mathbb{N}$,
$f^n$ is a conformal mapping from the component of
$f^{-n}(B_\infty(f^n(z_0),\delta))$ containing $z_0$ onto
$B_\infty(f^n(z_0),\delta)$, that is, the disk
$B_\infty(f^n(z_0),\delta)$ can be pulled back univalently along the
orbit of $z_0$. The radial (or conical) Julia set of $f(z)$, denoted
by $\mathcal{J}_f^r$, is defined to be the set of all conical
points. It is clear that a point $z_0\in\mathcal{J}_f\setminus
\mathcal{J}_f(\infty)$ is in $\mathcal{J}_f^r$ if
$\limsup\limits_{n\rightarrow\infty}d_\infty(f^n(z_0),\widehat{\mathcal{P}}(f))>0$.
Therefore, if $f(z)$ is in $\mathcal{H}(\widehat{\mathbb{C}})$, then
$\mathcal{J}_f^r=\mathcal{J}_f\setminus \mathcal{J}_f(\infty)$.

Set
\begin{equation}\label{3.2+}\mathcal{J}_f^{r,b}=\{z\in\mathcal{J}_f^r:\liminf_{n\rightarrow\infty}|f^n(z)|<+\infty\},\end{equation}
that is, $\mathcal{J}_f^{r,b}=\mathcal{J}_f\setminus
(\mathcal{I}(f)\cup\mathcal{J}_f(\infty)),$ and
$$\mathcal{J}_f^{r,u}=\mathcal{J}_f^r\cap \mathcal{I}(f),$$
where $$\mathcal{I}(f)=\{z\in \mathbb{C}: f^n(z)\rightarrow\infty\
{\rm as}\ n\rightarrow\infty\}.$$ Usually, $\mathcal{I}(f)$ is
called escaping set and it is also an important set which has
attracted many interests and researches. To some extent,
$\mathcal{I}(f)$ reflects the dynamic behavior of $f(z)$ on its
Julia set, especially, for transcendental entire functions. It is
well known that $\mathcal{J}_f=\partial\mathcal{I}(f)$. Generally,
we cannot have $\mathcal{I}(f)\subset\mathcal{J}_f$, while for a
function $f$ in Class $\mathcal{B}$, indeed the inclusion
$\mathcal{I}(f)\subset\mathcal{J}_f$ is true and therefore, for a
function $f\in\mathcal{H}(\widehat{\mathbb{C}})$,
$\mathcal{I}(f)\subset\mathcal{J}_f^r$ and
$\mathcal{J}_f^{r,u}=\mathcal{I}(f)$.

If $f(z)$ is in
$\mathcal{H}(\mathbb{C})\setminus\mathcal{H}(\widehat{\mathbb{C}})$,
then it is possible that $\mathcal{J}_f^{r,u}\not=\emptyset.$ For
that, the reader is referred to the proof of Dominguez
\cite{Dominguez} for existence of points in $\mathcal{I}(f)$ for a
transcendental meromorphic function with infinitely many poles. If
$f(z)$ is a transcendental meromorphic function in
$\mathcal{H}(\mathbb{C})$, then
\begin{equation}\label{1}\mathcal{J}_f^{r,b}
=\{z\in\mathcal{J}_f:\liminf_{n\rightarrow\infty}|f^n(z)|<+\infty\}.\end{equation}
If $f(z)$ is a transcendental entire function in
$\mathcal{H}(\mathbb{C})$, then
$\mathcal{J}_f^r=\mathcal{J}_f^{r,b}$. Consequently, the radial
Julia set is equivalently defined by the formula (\ref{1}) with
$\mathcal{J}_f^{r,b}$ replaced by $\mathcal{J}_f^r$, that is,
$\mathcal{J}_f^r=\mathcal{J}_f\setminus \mathcal{I}(f)$ for a
transcendental entire function in $\mathcal{H}(\mathbb{C})$.

It is easy to see that $\mathcal{J}_f^r\not=\emptyset$ for a
meromorphic function $f(z)$, as all repelling periodic points in the
Julia set must be in $\mathcal{J}_f^r$. In fact, as in the proof of
Theorem A of Stallard \cite{Stallard}, one has known that ${\rm
dim}_H(\mathcal{J}_f^{r,b})>0$. For a meromorphic function in
$\mathcal{H}(\mathbb{C})$, we can establish the following

\begin{thm}\label{thm3.3}\ Let $f(z)$ be in $\mathcal{H}(\mathbb{C})$. Then $${\rm
dim}_H(\mathcal{J}_f^r)={\rm dim}_H(\mathcal{J}_f^{r,b})=s(f)$$ and
furthermore, if $f\in\mathcal{H}(\widehat{\mathbb{C}})$, then ${\rm
dim}_H(\mathcal{J}_f)=s(f)$.\end{thm}

This result refines the result (2) in Theorem \ref{thm3.2} and
improves Theorem 2.7 of Kotus and Urbanski's \cite{KotusUrbanski}
for $f\in\mathcal{H}(\widehat{\mathbb{C}})$ in which they proved
that ${\rm dim}_H(\mathcal{J}_f)=s(f)$ if $f$ is so called strongly
regular(Note: We shall show that a function in
$\mathcal{H}(\widehat{\mathbb{C}})$ may not be strongly regular in
the remark (3) following Theorem \ref{thm3.5}). Theorem \ref{thm3.3}
follows from the following Lemma \ref{lem3.4} and Lemma
\ref{lem3.6}.

\begin{la}\label{lem3.4}\ Let $f(z)$ be in $\mathcal{H}(\mathbb{C}).$ Then ${\rm
dim}_H(\mathcal{J}_f^{r,b})\geq s(f)$.\end{la}

In order to prove Lemma \ref{lem3.4}, we need the following result
which is Lemma 4.1 of \cite{Stallard2}.

\begin{la}\label{lem3.5}\ Let $f(z)$ be a transcendental
meromorphic function in $\mathcal{H}(\mathbb{C}).$ Then for
$t<s(f)$, there exist some $a\in \mathcal{J}_f$ and $0<r<\delta$
with property that for arbitrarily large $C>0$, there exist
infinitely many values of $n\in\mathbb{N}$ such that
$$\sum_{f^n(z)=a,z\in B(a,r)}\frac{1}{(f^n)^\times(z)^t}>C.$$\end{la}

Now we begin to prove Lemma \ref{lem3.4}.

{\bf Proof of Lemma \ref{lem3.4}.}\ Let $t$ be a real number such
that $t<s(f)$. In view of Lemma \ref{lem3.5}, there exists a
$N\in\mathbb{N}$ such that
$$\sum_{f^N(z)=a,z\in B(a,r)}\frac{1}{|(f^N)'(z)|^t}\geq\left(\frac{1}{1+|a|^2}\right)^t
\sum_{f^N(z)=a,z\in B(a,r)}\frac{1}{(f^N)^\times(z)^t}>81^{2t}$$ and
$|(f^N)'(z)|>4\times 81$ for each $z\in B(a,r)$ with $f^N(z)=a$. For
such a $z$, we have an analytic branch $g$ of $f^{-N}$ on $B(a,4r)$
sending $a$ to $z$ and in view of the Koebe distortion theorem, for
$D=B(a,2r)$, we have
$${\rm diam}(g(D))<\frac{81{\rm diam}(D)}{
|(f^N)'(z)|}<\frac{1}{4}{\rm diam}(D)=r.$$ This implies that
$g(D)\subset D$. We write $\{D_i\}_1^M\ (1\leq M\leq\infty)$ for all
sets of $g(D)$ corresponding to the $z\in B(a,r)$ with $f^N(z)=a$
and by $g_i$ denote the branch of $f^{-N}$ such that $D_i=g_i(D)$.
Then $\{g_i\}_1^M$ is an iterated function system. Set
$$b_i=\inf\left\{\frac{|g_i(x)-g_i(y)|}{|x-y|}:\ x, y\in D\right\}.$$
In view of the Koebe distortion theorem, we can write
$$\frac{1}{81|(f^N)'(z_i)|}\leq b_i\leq \frac{81}{|(f^N)'(z_i)|}<\frac{1}{4},\ z_i=g_i(a).$$
Thus we have
$$\sum_{i=1}^Mb_i^t\geq \sum_{i=1}^M \frac{1}{(81|(f^N)'(z_i)|)^t}>1$$
and then we can choose some of $\{g_i\}_1^M$, say $\{g_i\}_1^P$ with
$1<P\leq M$ and $P<+\infty$ without any loss of generalities, and a
real number $t_0$ with $t\leq t_0$ such that
$$\sum_{i=1}^Pb_i^{t_0}=1.$$
In view of the well-known result (see Proposition 9.7 of
\cite{Falconer}), the iterated function system $\{g_i\}_1^P$
produces an invariant set $F$ for this system which has Hausdorff
dimension ${\rm dim}_H(F)\geq t_0\geq t$. It is obvious that
$F\subset \mathcal{J}_f^{r,b}$ and then ${\rm
dim}_H(\mathcal{J}_f^{r,b})\geq t$. Since $t$ is an arbitrary number
such that $t<s(f)$, we have ${\rm dim}_H(\mathcal{J}_f^{r,b})\geq
s(f)$. Lemma \ref{lem3.4} follows. \qed

\begin{la}\label{lem3.6}\ Let $f(z)$ be in $\mathcal{H}(\mathbb{C}).$ Then ${\rm
dim}_H(\mathcal{J}_f^{r,b})\leq s(f)$ and ${\rm
dim}_H(\mathcal{J}_f^{r,u})\leq s(f)$.\end{la}

{\bf Proof.}\ Take a sequence $\{a_m\}$ of complex numbers on
$\mathcal{J}_f$ so that $\mathcal{J}_f\subset \bigcup_{m=1}^\infty
B(a_m,\delta)$. Then
$$\mathcal{J}^{r,b}_f=\bigcup_{N=1}^\infty\bigcap_{M=1}^\infty\bigcup_{n=M}^\infty\bigcup_{m=1}^N
f^{-n}(B(a_m,\delta)).$$ We come to prove the equality. To the end,
we denote by $\mathcal{J}$ the set in the right side of the above
equality. Given arbitrarily a point $a\in \mathcal{J}^{r,b}_f$, we
have $f^{n_k}(a)\rightarrow c\not=\infty$ as $k\rightarrow\infty$
for a sequence of positive integers $\{n_k\}$. $c\in\mathcal{J}_f$
and for some $m_0$, $c\in B(a_{m_0},\delta)$. Obviously,
$a\in\bigcup_{m=1}^{m_0}f^{-n_k}(B(a_m,\delta))$ for large $n_k$ and
thus, for each $M\geq 1$,
$a\in\bigcup_{n=M}^\infty\bigcup_{m=1}^{m_0}f^{-n}(B(a_m,\delta))$.
This implies that $a\in\mathcal{J}$ and furthermore,
$\mathcal{J}^{r,b}_f\subseteq\mathcal{J}.$ Now given arbitrarily a
point $b\in\mathcal{J}$, there exists a $N\geq 1$ such that for all
$M\geq 1$, $b\in\bigcup_{n=M}^\infty\bigcup_{m=1}^N
f^{-n}(B(a_m,\delta))$ and for some $n_M\geq M$,
$b\in\bigcup_{m=1}^N f^{-n_M}(B(a_m,\delta))$, that is,
$f^{n_M}(b)\in\bigcup_{m=1}^NB(a_m,\delta)$. This implies that
$\liminf\limits_{n\rightarrow\infty}|f^n(b)|<+\infty$, that is,
$b\in\mathcal{J}_f^{r,b}$ and furthermore,
$\mathcal{J}\subseteq\mathcal{J}_f^{r,b}.$ We have attained our
desired result.

To observe the set $\mathcal{J}_f^{r,u}$, take a sequence
$\{\delta_m\}$ of positive numbers tending to zero. By $D_n^m$ we
denote the union of all components $V_n^m$ of
$f^{-n}(B_\infty(\infty,\delta_m))$ such that $f^n:V_n^m\rightarrow
B_\infty(\infty,\delta_m)$ is univalent. Set
$$\mathcal{I}_u=\bigcup_{N=1}^\infty\bigcap_{M=1}^\infty\bigcup_{m=N}^\infty\bigcup_{n=M}^\infty D_n^m.$$
Then $\mathcal{J}_f^{r,u}\subseteq\mathcal{I}_u$. Let us to prove
that. Fix arbitrarily a point $z_0\in\mathcal{J}_f^{r,u}$ and then
$f^n(z_0)\rightarrow\infty$ as $n\rightarrow\infty$ and we have a
$\delta_{m_0}$ and a sequence of positive integers $\{n_k\}$ such
that $f^{n_k}$ is a univalent analytic mapping from the component
$V_{n_k}^{m_0}$ containing $z_0$ onto
$B_\infty(\infty,\delta_{m_0})$. Hence $z_0\in\bigcup_{n=M}^\infty
D_n^{m_0}$ for each $M$ and furthermore, $z_0\in
\bigcup_{m=N}^\infty\bigcup_{n=M}^\infty D_n^m$ for each $M$ and
$N\leq m_0$. This easily implies that
$\mathcal{J}_f^{r,u}\subseteq\mathcal{I}_u$.

First of all, we come to prove that ${\rm
dim}_H(\mathcal{J}_f^{r,b})\leq s(f)$. Set
$$X_N=\bigcap_{M=1}^\infty\bigcup_{n=M}^\infty\bigcup_{m=1}^N
f^{-n}(B(a_m,\delta))$$ and thus
$\mathcal{J}_f^{r,b}=\bigcup_{N=1}^\infty X_N.$ It suffices to prove
that for each $N$, ${\rm dim}_H(X_N)\leq s(f)$ for our purpose.
Given arbitrarily $t>s(f)$, in view of the result (1) of Theorem
\ref{thm3.1} and the strictly decreasing property of $P(f,t)$ in
$t$, we have $P(f,t)<\lim\limits_{\tau\rightarrow
s(f)+0}P(f,\tau)\leq 0$ and there exists an integer $M_0>0$ such
that for $n\geq M_0$
$$\sum_{m=1}^N\mathcal{L}_t^n(1)(a_m)\leq
N\exp(nC)$$ for some fixed number $C$ with $P(f,t)<C<0$. Hence for
$M\geq M_0$
$$\sum_{n=M}^\infty\sum_{m=1}^N\mathcal{L}_t^n(1)(a_m)\leq
N\sum_{n=M}^\infty\exp(nC)=\frac{Ne^{MC}}{1-e^C}.$$

Let $g_{j,n,m}$ be an analytic branch of $f^{-n}$ on
$B(a_m,\delta)$. In view of the Koebe Theorem, we have
$${\rm diam}_\infty(g_{j,n,m}(B(a_m,\delta)))\leq K\delta
g^\times_{j,n,m}(a_m)=\frac{K\delta}{(f^n)^\times(g_{j,n,m}(a_m))}.$$
This implies that
\begin{eqnarray*}\sum_{n=M}^\infty\sum_{m=1}^N\sum_{j=1}^\infty{\rm
diam}_\infty(g_{j,n,m}(B(a_m,\delta))^t&\leq&
(K\delta)^t\sum_{n=M}^\infty\sum_{m=1}^N\sum_{f^n(z)=a_m}\frac{1}{((f^n)^\times(z))^t}\\
&\leq& (K\delta)^t\frac{Ne^{MC}}{1-e^C}.\end{eqnarray*}

Since $g_{j,n,m}(\zeta)$ uniformly converges on $B(a_m,\delta)$ to a
point on $\widehat{\mathbb{C}}$, it is easy to see that $${\rm
diam}_\infty(g_{j,n,m}(B(a_m,\delta)))\rightarrow 0\ (
n\rightarrow\infty).$$ By noting that
$X_N\subset\bigcup_{n=M}^\infty\bigcup_{m=1}^N
f^{-n}(B(a_m,\delta))$ for each $M>0$, we have the Hausdorff measure
$\mathcal{H}^t(X_N)=0$. This yields $t\geq {\rm dim}_H(X_N)$ and so
${\rm dim}_H(X_N)\leq s(f)$. This is our desired result.

Write $z_{n,m}$ for the point in $V_n^m$ such that
$f^n(z_{n,m})=\infty$. It is easy to prove that for $t>s(f)$,
$$\lim_{n\rightarrow\infty}\sum_{z_{n,m}}\frac{1}{
(f^n)^\times(z_{n,m})^t}=0$$ and
$P_\infty(f,t)=\limsup\limits_{n\rightarrow\infty}\frac{1}{n}\log
\sum\limits_{z_{n,m}}\frac{1}{(f^n)^\times(z_{n,m})^t}$ is strictly
decreasing in $t$. The same argument as above can be used to prove
that ${\rm dim}_H(\mathcal{J}_f^{r,u})\leq {\rm
dim}_H(\mathcal{I}_u)\leq s(f)$. Thus we complete the proof of Lemma
\ref{lem3.6}.\qed

Therefore, from Theorem \ref{thm3.3} it follows that ${\rm
dim}_H(\mathcal{J}_f)=s(f)$ if and only if ${\rm
dim}_H(\mathcal{I}(f))\break\leq s(f)$. Urbanski and Zdunik
\cite{UrbanskiZdunik} proved that ${\rm dim}_H(\mathcal{J}_f^r)<2$
for $f(z)=\lambda e^z$ with $\lambda\in\mathbb{C}\setminus\{0\}$
such that $f(z)$ has an attracting periodic point. Obviously, such a
function is hyperbolic on the complex plane and in fact, in view of
Proposition 2.1 of Urbanski and Zdunik \cite{UrbanskiZdunik1}, it is
basically hyperbolic with respect to the Euclidean metric. On the
other hand, McMullen \cite{McMullen1} proved that ${\rm
dim}_H(\mathcal{I}(f))=2$. Therefore, even if a function is
hyperbolic with respect to the Euclidean metric, its Hausdorff
dimension may be large than $s(f)$. The case ${\rm
dim}_H(\mathcal{I}(f))<{\rm dim}_H(\mathcal{J}_f^r)$ is also
possible. Consider the function $g(z)=\lambda \tan z$. In Example 6
of \cite{KotusUrbanski1}, Kotus and Urbanski proved that ${\rm
dim}_H(\mathcal{I}(g))\leq\frac{1}{2}<{\rm dim}_H(\mathcal{J}_g^r)$.
For suitable $\lambda$, $g\in\mathcal{H}(\widehat{\mathbb{C}})$.
Observing the examples discussed in \cite{KotusUrbanski1}, we
propose a question: if $\infty$ is an asymptotic value of a
meromorphic function $f(z)$, should we always have that ${\rm
dim}_H(\mathcal{I}(f))\geq{\rm dim}_H(\mathcal{J}_f^r)$? Or if
$\infty$ is not an asymptotic value of a meromorphic function
$f(z)$, should we always have that ${\rm
dim}_H(\mathcal{I}(f))\leq{\rm dim}_H(\mathcal{J}_f^r)$?

Next in view of results of Walters \cite{Walters} as Kotus and
Urbanski \cite{KotusUrbanski} did (compare \cite{Zheng1}), we
consider thermodynamical formalism of hyperbolic meromorphic
functions on the Riemann sphere. The reader is referred to Mayer and
Urbanski \cite{MayerUrbanski} for thermodynamical formalism of some
of meromorphic functions of finite order hyperbolic with respect to
the Euclidean metric, but not on the Riemann sphere.

We first of all are concerned for what $t$ the operator
$\mathcal{L}_t$ is a mapping of
$C(\widehat{\mathcal{J}}_f)\rightarrow C(\widehat{\mathcal{J}}_f)$.
Lemma 7.3.3 of \cite{Zheng1} actually asserts the existence of a
real number $t(f)\leq s(f)$ such that $\mathcal{L}_t$ is a mapping
of $C(\widehat{\mathcal{J}}_f)\rightarrow
C(\widehat{\mathcal{J}}_f)$ only when $t(f)<t$. The following result
confirms the case when $t=s(f)$. When $\mathcal{L}_t:
C(\widehat{\mathcal{J}}_f)\rightarrow C(\widehat{\mathcal{J}}_f)$,
we denote by $\mathcal{L}^*_t$ the dual operator of $\mathcal{L}_t$.

\begin{la}\label{lem3.7}\ Let $f(z)$ be a transcendental
meromorphic function in $\mathcal{H}(\widehat{\mathbb{C}}).$ Then
for $s=s(f)$, we have
\begin{equation}\label{3.11}T=\sup\{\sum_{f(z)=a}\frac{1}{f^\times(z)^s}: a\in
\widehat{\mathcal{J}}_f\}<+\infty,\end{equation} $P(f,s)=0$ and
furthermore, $\mathcal{L}_s$ is a continuous operator of
$C(\widehat{\mathcal{J}}_f)\rightarrow
C(\widehat{\mathcal{J}}_f)$.\end{la}
%Then $\forall \varepsilon>0$, $\exists \eta$ with $0<\eta<\delta$
%such that for $a, b\in\mathcal{J}_f$, when $d(a,b)<\eta$, we have
%\begin{equation}\label{3.10}\sum_{f(z)=a}\left|\frac{1}{f^\times(z)^t}-\frac{1}{f^\times(z')^t}\right|<\varepsilon\end{equation}
%where $f(z')=b$ such that $z$ and $z'$ are in a component of
%$f^{-1}(B(a,\eta))$.

{\bf Proof.}\ To prove (\ref{3.11}). Take $M$ points $a_m\in
\mathcal{J}_f\ (m=1,2,...,M)$ such that
$\widehat{\mathcal{J}}_f\subset \bigcup_{m=1}^MB(a_m,\delta)\cup
B_\infty(\infty,\delta)$ and take a positive integer $N$ such that
for an arbitrary pair of $m$ and $n$, $B(a_m,\delta)\cap
f^{-N+1}(a_n)\not=\emptyset$, and $B_\infty(\infty,\delta)\cap
 f^{-N+1}(a_n)\not=\emptyset$, where $a_n=\infty$ for $n=M+1$. For
$t>s(f)$, in view of (\ref{3.2}), we have a positive integer
$P=P(t)$ such that
\begin{eqnarray}\label{3.12}1&>&\sum_{f^{PN}(z)=a}\frac{1}{
((f^{PN})^\times(z))^t}\nonumber\\
&=&\sum_{f^{N}(w)=a}\frac{1}{
((f^{N})^\times(w))^t}\sum_{f^{(P-1)N}(z)=w}\frac{1}{
((f^{(P-1)N})^\times(z))^t}.\end{eqnarray} If for each $w$,
$\sum_{f^{(P-1)N}(z)=w}\frac{1}{((f^{(P-1)N})^\times(z))^t}\geq 1,$
then we have $\sum_{f^{N}(w)=a}\frac{1}{((f^{N})^\times(w))^t}<1$;
If for some $w$, $\sum_{f^{(P-1)N}(z)=w}\frac{1}{
((f^{(P-1)N})^\times(z))^t}<1,$ then in the same step as above, we
shall have at last $\sum_{f^{N}(w)=b}\frac{1}{
((f^{N})^\times(w))^t}<1$ for some $b=b(t)\in\mathcal{J}_f$. There
exists an integer $Q$ with $1\leq Q\leq M$ such that $b\in
B(a_Q,\delta)$ or $b\in B_\infty(\infty,\delta)$. Assume without any
loss of generalities that $b\in B(a_Q,\delta)$. In view of the Koebe
distortion theorem, we have
$$\sum_{f^{N}(w)=a_Q}\frac{1}{
((f^{N})^\times(w))^t}\leq\sum_{f^{N}(w)=b}\frac{K^t}{
((f^{N})^\times(w))^t}<K^t$$ for an absolute constant $K>1$ and
furthermore, it follows from (\ref{3.9}) that
$$\mathcal{L}_t(1)(w)\leq \mathcal{L}_t^N(1)(a_Q)(f^{N-1})^\times(w)^t
\leq K^t(f^{N-1})^\times(w)^t$$ for all $w\in f^{-N+1}(a_Q)$. Take a
point $w_n^m\in B(a_m,\delta)\cap f^{-N+1}(a_n)$ for each pair of
$m$ and $n$ and a point $w_n^{M+1}\in B_\infty(\infty,\delta)\cap
f^{-N+1}(a_n)$ where $a_n=\infty$ for $n=M+1$.

Set
$$T(t)=\sup\{K^{2t}(f^{N-1})^\times(w_n^m)^t:\ m, n=1,2,...,M+1\}.$$
For an arbitrary point $c\in\widehat{\mathcal{J}}_f$, we have $c\in
B(a_m,\delta)$ for some $m$ or $c\in B_\infty(\infty,\delta)$ and in
view of the Koebe distortion theorem, we have
$$\mathcal{L}_t(1)(c)\leq
K^t\max\{\mathcal{L}_t(1)(w_Q^m), \mathcal{L}_t(1)(w_Q^{M+1})\}\leq
T(t).$$ Letting $t\rightarrow s(f)+0$ implies that
\begin{equation}\label{3.13}\mathcal{L}_s(1)(c)\leq
T(s),\end{equation} and then we have proved that (\ref{3.13}) holds
uniformly on $\widehat{\mathcal{J}}_f$, that is, (\ref{3.11}) holds.

Using the result (iv) in Lemma \ref{lem3.1+} to (\ref{3.13}) yields
$P(f,s)\leq T(s)$ and it is easy to see that $P(f,s)=0$.

In view of Lemma \ref{lem3.3+}, we have a positive number
$\eta_0<\delta$ such that for each $a\in \mathcal{J}_f$, $f^{-1}$
can be divided into single-valued analytic branches on $B(a,\delta)$
and for each point $b\in \mathcal{J}_f$ with $d(a,b)<\eta_0$, we
have
\begin{eqnarray*}\sum_{f(z)=a}\left|\frac{1}{f^\times(z)^s}-\frac{1}{f^\times(z')^s}\right|&=&
\sum_{f(z)=a}\frac{1}{f^\times(z)^s}\left|1-\frac{f^\times(z)^s}{f^\times(z')^s}\right|\\
&\leq&\sum_{f(z)=a}\frac{1}{f^\times(z)^s}K_sd(a,b)\\
&\leq& T(s)K_sd(a,b),\end{eqnarray*} where $z'=f^{-1}_z(b)$ and
$f^{-1}_z$ is the branch of $f^{-1}$ sending $a$ to $z$; When
$a=\infty$, we also have the above inequality with all $d$ replaced
by $d_\infty$, that is, we use the sphere metric.

For each $g\in C(\widehat{\mathcal{J}}_f)$, noticing that
$\widehat{\mathcal{J}}_f$ is compact and $||g||=\sup\{|g(z)|:z\in
\widehat{\mathcal{J}}_f\}<+\infty$, for an arbitrary pair of two
points $a,b\in\mathcal{J}_f$ with $d(a,b)<\eta_0$, we have
\begin{eqnarray*}
|\mathcal{L}_sg(a)-\mathcal{L}_sg(b)|&\leq&\sum_{f(z)=a}\left|\frac{g(z)}{
f^\times(z)^s}-\frac{g(z')}{f^\times(z')^s}\right|\\
&\leq&\sum_{f(z)=a}\frac{1}{f^\times(z)^s}|g(z)-g(z')|\\
&+&|g(z')|\sum_{f(z)=a}\left|\frac{1}{f^\times(z)^s}-\frac{1}{
f^\times(z')^s}\right|\\
&\leq&
T(s)\sup_{f(z)=a}|g(z)-g(z')|+||g||T(s)K_sd(a,b);\end{eqnarray*}
When $a=\infty$, we also have the above inequality with $d$ replaced
by $d_\infty$. In view of the Koebe distortion theorem, we have
\begin{eqnarray*}
d_\infty(z,z')&=&d_\infty(f^{-1}_z(a),f^{-1}_z(b))\leq
K(f^{-1})^\times(a)d_\infty(a,b)\\
&=&\frac{K}{f^\times(z)}d_\infty(a,b)\leq
CT(s)^{1/s}d_\infty(a,b).\end{eqnarray*} Therefore,
$\mathcal{L}_sg(w)$ is continuous on $\widehat{\mathcal{J}}_f$.

Thus we complete the proof of Lemma \ref{lem3.7}. \qed

In view of Lemma \ref{lem3.7}, the existence of $\mathcal{L}_s$
together with $P(f,s)=0$ implies the existence of an unique
conformal measure $\mu_s$ for a
$f\in\mathcal{H}(\widehat{\mathbb{C}})$, which is the fixed point of
$\mathcal{L}^*_s$, and then in terms of the developed results of
Walters \cite{Walters} in \cite{Zheng2}, i.e., Theorem 2.6 of
\cite{Zheng2}, applying the expanding property stated in Theorem
\ref{thm2.3}, Lemma \ref{lem3.3+} and the summable property of
$\mathcal{L}_s(1)(w)$ yields the existence of an unique invariant
measure $m_s$ which is equivalent to $\mu_s$. In one word, we can
establish the following

\begin{thm}\label{thm3.5}\ \ Let $f(z)$ be a transcendental
meromorphic function in $\mathcal{H}(\widehat{\mathbb{C}}).$ Then

(1) $s=s(f)={\rm dim}_H(\widehat{\mathcal{J}}_f)<2$;

(2) there exists a $f^\times(z)^s$-conformal measure $\mu_s$ of
$f(z)$ on $\mathcal{J}_f\setminus f^{-1}(\infty)$ with
$\mu_s(\mathcal{J}_f\setminus f^{-1}(\infty))=1$ and $\mu_s$ is
positively nonsingular and nonsingular, that is, $\mu_s\circ
f\ll\mu_s$ and $\mu_s\circ f^{-1}\ll\mu_s$;

(3) there exists a $h_s\in\mathcal{C}(\widehat{\mathcal{J}}_f)$ with
$h_s>0$ such that $\mu_s(h_s)=1$ and $\mathcal{L}_s(h_s)=h_s$ and
$m_s=h_s\mu_s$ is an invariant Gibbs measure for $f(z)$;

(4) both of $\mu_s$ and $m_s$ are positive on any open subset of
$\mathcal{J}_f$ and atomless;

(5) for every $g\in \mathcal{C}(\mathcal{J}_f)$, we have
$$E_{\mu_s}(g|f^{-1}\varepsilon)(x)=\mathcal{L}_s(g)\circ f(x),\
\mu_s-{\rm a.e},$$ and for every invariant measure $\sigma$ of
$f(z)$, we have
$$0=\mu_s(I_{\mu_s}(\varepsilon|f^{-1}\varepsilon)+\varphi_s)\geq
\sigma(I_{\sigma}(\varepsilon|f^{-1}\varepsilon)+\varphi_s);$$

(6) if $\log f^\times (z)$ is $\mu_s$-integrable, then
$$s=\lim_{r\rightarrow 0}\frac{\log\mu_s(B_\infty(a,r))}{\log r}, {\rm
for}\ a\in \mathcal{J}_f, \mu_s-{\rm a.e}$$ and
$$s=s(f)={\rm dim}_H(\widehat{\mathcal{J}}_f)={\rm dim}_P(\widehat{\mathcal{J}}_f),$$
 and $0<\chi_{\mu_s}(f)<\infty$ and furthermore
$$s=s(f)=\frac{H_{\mu_s}(\varepsilon|f^{-1}\varepsilon)}{\chi_{\mu_s}(f)};$$

(7) $\mathcal{H}^s(\mathcal{J}_f)<+\infty$,
$\mathcal{P}^s(\mathcal{J}_f)>0$, $\mathcal{H}^s\ll\mu_s$ and
$\frac{\mathrm{d}\mathcal{H}^s}{\mathrm{d}\mu_s}<+\infty,$ where
$\mathcal{H}^s$ is the Hausdorff measure and $\mathcal{P}^s$ the
packing measure of $s$ dimension;

(8) $(f,\mu_s)$ is an exact endomorphism.
\end{thm}

We remark on Theorem \ref{thm3.5}. (1) The result (1) is independent
of the results of Walters \cite{Walters}, and the equality
"$s(f)={\rm dim}_H(\widehat{\mathcal{J}}_f)$" is proved in Theorem
2.7 in \cite{KotusUrbanski} under additional assumption of that $f$
is strongly regular.

(2) In \cite{KotusUrbanski}, Kotus and Urbanski stated, in fact, the
existence of conformal measure and invariant measure of
transcendental meromorphic functions in
$\mathcal{H}(\widehat{\mathbb{C}})$, which was extracted from the
results of Walters \cite{Walters}, by assuming that $\varphi_s$ is
summable and $P(f,s)=0$. However, in Lemma \ref{lem3.7}, we proved
that $\varphi_s$ is summable over $\mathcal{J}_f$, $P(f,s)=0$ and
$\mathcal{L}_s$ can be extended to a linear operator  from
$\mathcal{C}(\widehat{\mathcal{J}}_f)$ to itself, which confirms the
existence of $\mu_s$ and in terms of Theroem \ref{thm2.3}, $f^N$ is
expanding over $\mathcal{J}_f\setminus \mathcal{J}_f(\infty)$ for a
fixed $N$, thus Theorem 2.6 of \cite{Zheng2} yields the existence of
$m_s$ and other results.

(3) It was proved in the Claim in the proof of Theorem 2.7 in
\cite{KotusUrbanski} that if $s>\tau(f)=\inf\{t\geq 0:
P(f,t)<\infty\}$, namely $f$ is strongly regular, then $\log
f^\times(z)$ is $\mu_s$-integrable and $0<\chi_{\mu_s}(f)<\infty$ so
that $s(f)={\rm dim}_H(\widehat{\mathcal{J}}_f).$ However, there
exist meromorphic functions in $\mathcal{H}(\widehat{\mathbb{C}})$
such that $\log f^\times(z)$ is not $\mu_s$-integrable. Therefore,
the result (1) is a generalization of Theorem 2.7 in
\cite{KotusUrbanski}. Actually, there exist functions in
$\mathcal{H}(\widehat{\mathbb{C}})$ such that the Hausdorff
dimensions of their Julia sets do not equal to their packing
dimensions. Put
$$g(z)=\frac{\tan z}{(2m)^m\prod_{j=1}^m(z-j\pi)},\ {\rm where}\
m\in\mathbb{N}.$$ Stallard \cite{Stallard1} proved that for $m>8$,
${\rm dim}_P(\mathcal{J}_g)\geq\frac{1}{2}>\frac{4}{m}\geq{\rm
dim}_H(\mathcal{J}_g)$ and $\mathcal{P}(g)\cap
\widehat{\mathcal{J}}_g=\emptyset$. It is obvious that $\infty$ is
not a singular value of $g(z)$ and so $g(z)$ is in
$\mathcal{H}(\widehat{\mathbb{C}})$. In view of Theorem
\ref{thm3.5}, $\log g^\times(z)$ is not $\mu_s$-integrable on
$\mathcal{J}_g$ and $s(g)=\tau(g)$.

(4) For the function $f_{p,\lambda}$ in (\ref{1}), we have
$s(f_{p,\lambda})>t(f_{p,\lambda})$, and then
$\log(f_{p,\lambda})^\times$ is $\mu_s$-integrable on
$\mathcal{J}_{f_{p,\lambda}}.$

(5) Since $P(f,s)=0$ is proved, all results in \cite{KotusUrbanski}
for the regular Walters expanding conformal map about the Hausdorff
measure and packing measure apply to the functions in
$\mathcal{H}(\widehat{\mathbb{C}})$. The result (7) of Theorem
\ref{thm3.5} follows from Theorem 2.13 in \cite{KotusUrbanski}
together with $P(f,s)=0$.

\vskip 0.5cm {\bf Acknowledgement}\ \ The author is grateful to
Walter Bergweiler for his value discussion and hospitality. And he
also wishes to thank the European Science Foundation for the
support.

\end{document}